\documentclass[11pt]{amsart}                  
\usepackage{amssymb,amsmath,amsfonts,latexsym}
\usepackage[margin=1.1in]{geometry}
\usepackage{enumerate}

\newcommand{\reals}{{\mathbb R}}

\newcommand{\cala}{{\mathcal A}}

\newcommand{\calh}{{\mathcal H}}

\newcommand{\calw}{{\mathcal W}}

\theoremstyle{plain}
        \newtheorem{theorem}{Theorem}[section]
        \newtheorem{lemma}[theorem]{Lemma}
        \newtheorem{remark}[theorem]{Remark}
        \newtheorem{proposition}[theorem]{Proposition}

\begin{document}
\title[universal deformation formula]{A universal deformation
formula for Connes-Moscovici's Hopf algebra without any projective
structure}
\author{Xiang Tang and Yijun Yao}
\maketitle
\begin{abstract}
We construct a universal deformation formula for Connes-Moscovici's
Hopf algebra without any projective structure using Fedosov's
quantization of symplectic diffeomorphisms.
\end{abstract}
\section{Introduction}
In the study of index theory of a transverse elliptic differential
operator of a codimension one foliation, Connes and Moscovici
discovered a Hopf algebra $\calh_1$ which governs the local
symmetry in computing the Chern character. In this paper, we study
deformation theory of this Hopf algebra. In particular, we prove
that the Hopf algebra $\calh_1$ has a universal deformation
formula.

In \cite{CM03-2}, inspired from Rankin-Cohen brackets on modular
forms, Connes and Moscovici constructed a universal deformation
formula of the Hopf algebra $\calh_1$ with a projective structure.
By a universal deformation formula of a Hopf algebra $A$, we mean
an element $R\in A[[\hbar]]\otimes A[[\hbar]]$ satisfying
\[
(\Delta\otimes 1)R (R\otimes 1)=(1\otimes \Delta)R (1\otimes R),\
\text{and}\ \epsilon\otimes 1(R)=1\otimes 1=1\otimes \epsilon (R).
\]
In \cite{bty}, we together with Bieliavsky provided a geometric
interpretation of a projective structure in the case of a
codimension one foliation. And as a result, we (with Bieliavsky)
obtained a geometric way to reconstruct Connes-Moscovici's
universal deformation formula. Furthermore, a new and interesting
result we proved in \cite{bty}[Prop. 6.1] is that even without a
projective structure, the first Rankin-Cohen bracket
\[
RC_1=S(X)\otimes Y+Y\otimes X\in \calh_1\otimes \calh_1
\]
is a noncommutative Poisson structure, i.e. $RC_1$ is a Hochschild
cocycle and $(1\otimes \Delta)RC_1(1\otimes RC_1)-(\Delta\otimes
1)RC_1(RC_1\otimes 1)$ is a Hochschild coboundary. This leads to a
question whether $\calh_1$ has a universal deformation formula
without a projective structure.

In this article, we give a positive answer to the above question
and introduce a geometric construction of a universal deformation
formula of $\calh_1$ without a projective structure. The idea of
this construction goes back to Fedosov \cite{fe} in his study of
deformation quantization of a symplectic diffeomorphism. Fedosov
developed in \cite{fe} a systematic way to quantize a symplectic
diffeomorphism to an endomorphism of the quantum algebra no matter
whether it preserves or not the chosen symplectic connection.
Fedosov also observed that composition of quantized symplectic
diffeomorphisms does not preserve associativity. Instead, it
satisfies a weaker associativity property---associative up to an
inner endomorphism. This picture can be explained using the
language of ``gerbes and stacks" as \cite{bgnt}. In any case,
Fedosov's construction does give rise to a deformation
quantization of the groupoid algebra associated to a pseudogroup
action on a symplectic manifold.

In this paper, we apply this idea to the special case that the
symplectic manifold is $\reals\times \reals^+$ and the Poisson
structure is $\partial_x\wedge \partial_y$, with $x$ the
coordinate on $\reals$ and $y$ the coordinate on $\reals^+$. We
consider symplectic diffeomorphisms on $\reals\times \reals^+$ of
the form
\[
\gamma:(x,y)\to \left(\gamma(x),\frac{y}{\gamma'(x)}\right),
\]
where $\gamma$ is a local diffeomorphism on $\reals$. In this
case, Fedosov's construction of quantization of symplectic
diffeomorphism can be computed explicitly. In particular, we are
able to prove that the resulting star product on the groupoid
algebra $C_c^\infty(\reals\times\reals^{+})\rtimes \Gamma$ can be
expressed by
\[
f\alpha\star g\beta=m(R(f\alpha\otimes g\beta)),
\]
where $m$ is the multiplication map on $C_c^\infty(\reals\times
\reals^+)\rtimes \Gamma$, and $R$ is an element in
$\calh_1[[\hbar]]\otimes \calh_1[[\hbar]]$. An important property
is that the $\calh_1$ action on the collection of all
$C_c^\infty(\reals\times \reals^+)\rtimes \Gamma$ for all
pseudogroup $\Gamma$ is fully injective because this action is
equivalent to the action used by Connes and Moscovici to define
$\calh_1$. With this observation, we can derive all the property
of $R$ as a universal deformation formula from the corresponding
properties about the star product on $C_c^\infty(\reals\times
\reals^+)\rtimes \Gamma$.

The notion of universal deformation formula of a Hopf algebra is
closely related to the solution of quantum Yang-Baxter equation.
We hope that our construction will shed a light on the study of
deformation theory of the Hopf algebra $\calh_1$ and also
codimension one foliations.

This article is organized as follows. We review in section 2
Fedosov's theory of deformation quantization of symplectic
diffeomorphisms. We provide a detailed proof of the fact that this
defines a deformation of the groupoid algebra
$C_c^\infty(\reals\times \reals^+)\rtimes \Gamma$. In section 3, we
prove the main theorem of this paper that $\calh_1$ has a universal
deformation formula using Fedosov's theory reviewed in section 2. In
section 4, we compute explicitly our universal deformation formula
up to $\hbar^2$. We observe that when the $\calh_1$ action is
projective, the universal deformation formula obtained in this paper
does not agree with the one introduced by Connes and Moscovici in
\cite{CM03-2}. Instead, our lower order term computation suggests
that in the case of a projective action, these two universal
deformation formulas should be related by an isomorphism expressed
by elements in $\calh_1[[\hbar]]$ and the projective structure
$\Omega$. In the appendix we discuss the associativity of the
Eholzer product on modular forms, which was used by Connes and
Moscovici in constructing their Rankin-Cohen deformation.
\\

\noindent{\bf Acknowledgements:} We would like to thank A. Connes
and H. Moscovici for explaining the Hopf algebra $\calh_1$ and
Rankin-Cohen brackets. The first author would like to thank A.
Gorokhovsky, R. Nest, and B. Tsygan for explaining Fedosov's
quantization of symplectic diffeomorphisms and their ideas of
deformation of groupoid algebras. The research of the first author
is partially supported by NSF Grant 0703775.
\section{Quantization of symplectic diffeomorphisms}
In this section, we briefly recall Fedosov's construction of
quantization of a symplectic diffeomorphism. And we use this idea
to define a deformation of a groupoid algebra coming from a
pseudogroup action on a symplectic manifold. We learned this
construction from A. Gorokhovsky, R. Nest, and B. Tsygan.

In Fedosov's approach to deformation quantization of A symplectic
manifold $(M, \omega)$, a flat connection $D$ (also called Fedosov
connection) on the Weyl algebra bundle $\calw$ plays an essential
role. The quantum algebra is identified with the space of flat
sections $\calw_D$ of $\calw$. A question arises when one wants to
quantize a symplectic diffeomorphism. Because a symplectic
diffeomorphism may not preserve $D$, the canonical lifting of a
symplectic diffeomorphism to the Weyl algebra bundle $\calw$ may
not act on the quantum algebra $\calw_D$. How can we quantize a
symplectic diffeomorphism in this case? Fedosov in \cite{fe}
studied this question. The answer he came up with fits well with
the language of ``stack of algebras". In the following we briefly
review Fedosov's results \cite{fe}[Section 4].

A symplectic diffeomorphism $\gamma:M\to M$ naturally acts on the
cotangent bundle $\gamma:T^*M\to T^*M$. Therefore $\gamma$ lifts
to an endomorphism on the Weyl algebra bundle $\gamma:\calw\to
\calw$. It is easy to check that if $\gamma(D):=\gamma \circ
D\circ \gamma^{-1}=D$, then $\gamma$ defines an algebra
endomorphism on the quantum algebra $\calw_D=\ker(D)$, which is
called a quantization of the symplectic diffeomorphism $\gamma$.
We with Bieliavsky in \cite{bty} used this idea to construct a
universal deformation formula of $\calh_1$ with a projective
structure.

The quantization of $\gamma$ when $\gamma(D)\ne D$ is more
involved. Fedosov in \cite{fe}[Sec. 4] purposed the following
construction of quantization. We start with extending the standard
Weyl algebra $W$ to $W^+$,
\begin{enumerate}
\item an element $u$ of $W^+$ can be written as
\[
u=\sum_{2k+l\geq0}\hbar^ka_{k,i_1,\cdots, i_l}y^{i_1}\cdots
y^{i_l}
\]
where $(y^1, \cdots, y^{2n})$ are coordinates on the standard
symplectic vector space $(V, \omega)$. In the above sum, we allow
$k$ to be negative.
\item There are a finite number of terms with a given
total degree $2k+l\geq0$.
\end{enumerate}
We remark that the Moyal product extends to a well define product
on $W^+$. And we consider the corresponding extension $\calw^+$ of
the Weyl algebra bundle $\calw$ associated to $W^+$. Write
$\gamma(D)=D+i/\hbar[\Delta_\gamma, \cdot]$ with $\Delta_\gamma$ a
smooth section of the Weyl algebra bundle $\calw$. We remark that
in principle, we can change $\Delta_\gamma$ by any section in the
center of $\calw$. But when we fix the expression of
$D=d+i/\hbar[r,\cdot]$, $\Delta_\gamma$ has a canonical choice
$\gamma(r)-r$. In the following of this paper, we will always work
with this choice of $\Delta_\gamma$. We consider the following
equation
\begin{equation}\label{eq:u}
DU_\gamma=-\frac{i}{\hbar}\Delta_\gamma\circ U_\gamma,
\end{equation}
where $U_\gamma$ is an invertible section of $\calw^+$. Fedosov
\cite{fe}[Thm 4.3] proved that equation (\ref{eq:u}) always has
solutions. In general, solutions to equation (\ref{eq:u}) are not
unique. But the following induction procedure
\[
U_{n+1}=1+\delta^{-1}\{(D+\delta)U_n+(i/\hbar)\Delta_\gamma\circ
U_n\},\ \ \ \ \ U_0=1
\]
uniquely determines an invertible solution to Equation
(\ref{eq:u}). By this induction, we see that $U$ is a solution to
the following equation
\begin{equation}\label{eq:induction-u}
U=1+\delta^{-1}\{(D+\delta)U+i/\hbar \Delta_\gamma\circ U\},
\end{equation}
which actually has a unique solution because
$\delta^{-1}\{(D+\delta)U+i/\hbar \Delta_\gamma\circ U\}$ raises
the total degree of $U$ by 1. We will always work with this
solution in this paper. By Equation (\ref{eq:u}), $U_\alpha^{-1}$
satisfies the following equation,
\[
\begin{split}
DU_\alpha ^{-1}&=-U^{-1}_{\alpha}\circ DU_\alpha \circ
U^{-1}_\alpha\\
&=\frac{i}{\hbar}U_\alpha^{-1}\circ \Delta_\alpha.
\end{split}
\]
Using Equation (\ref{eq:induction-u}), it is not difficult to
check that $U_\alpha^{-1}$ is the unique solution to the following
equation,
\begin{equation}\label{eq:u-inverse}
V=1+\delta^{-1}\{(D+\delta)V-i/\hbar V\circ \Delta_\alpha\},
\end{equation}
which can be constructed by same induction as above.

We have the following lemma on $\Delta_\alpha$.
\begin{lemma}
\label{lem:delta}
\begin{enumerate}
\item
$\Delta_{\alpha}+\alpha(\Delta_{\beta})=\Delta_{\alpha\beta}$,
\item $\alpha(\Delta_{\alpha^{-1}})=-\Delta_\alpha$.
\end{enumerate}
\end{lemma}
\begin{proof}
(1) $D+i/\hbar[\Delta_{\alpha\beta},
\cdot]=\alpha\beta(D)=\alpha(\beta(D))=\alpha(D+i/\hbar[\Delta_\beta,\cdot])
=\alpha(D)+i/\hbar[\alpha(\Delta_\beta),\cdot]=D+i/\hbar[\Delta_\alpha,\cdot]
+i/\hbar[\alpha(\Delta_\beta),\cdot]=D+i/\hbar[\Delta_\alpha+\alpha(\Delta_\beta),\cdot]$.
Therefore by the defining property of $\Delta_{\alpha\beta}$ and
its uniqueness, we conclude that
$\Delta_{\alpha\beta}=\Delta_\alpha+\alpha(\Delta_\beta)$.

(2) Corollary of (1).
\end{proof}

We will need the following properties of $U_\alpha$ later in our
construction.
\begin{proposition}
\label{prop:u} The assignment $\alpha\mapsto U_\alpha$ satisfies
the following properties,
\begin{enumerate}
\item $D(\alpha(U_\beta))=-i/\hbar\big(\Delta_{\alpha\beta}\big)\circ \alpha(U_{\beta})+
i/\hbar \alpha(U_\beta)\circ\Delta_{\alpha}$;
\item $\alpha(U_{\alpha^{-1}})=U_{\alpha}^{-1}$.
\end{enumerate}
\end{proposition}
\begin{proof}
(1)We can use Equation (\ref{eq:induction-u}) to prove a stronger
statement. We compute
\begin{eqnarray*}
\alpha({U_\beta})&=&\alpha(1+\delta^{-1}\{(D+\delta){U_\beta}+i/\hbar\Delta_\beta\circ
{U_\beta}\})\\
&=&1+\delta^{-1}\{(\alpha(D)+\delta)\alpha({U_\beta})+i/\hbar\alpha(\Delta_\beta)\circ
\alpha({U_\beta})\}\\
&=&1+\delta^{-1}\{(D+\delta)\alpha({U_\beta})+i/\hbar[\Delta_\alpha,\alpha({U_\beta})]
+i/\hbar\alpha(\Delta_\beta)\circ\alpha({U_\beta})\}\\
&=&1+\delta^{-1}\{(D+\delta)\alpha({U_\beta})+i/\hbar(\Delta_\alpha+\alpha(\Delta_\beta))\circ
\alpha({U_\beta})-i/\hbar \alpha({U_\beta})\circ \Delta_\alpha\}.
\end{eqnarray*}
By applying Lemma \ref{lem:delta} to the last line, we conclude
that $\alpha(U_\beta)$ is the unique solution to the following
equation
\begin{equation}\label{eq:alpha-U-beta}
\alpha(U_\beta)=1+\delta^{-1}\{(D+\delta)\alpha(U_\beta)+i/\hbar\Delta_{\alpha\beta}\circ
\alpha(U_\beta)-i/\hbar \alpha(U_\beta)\circ \Delta_\alpha\}.
\end{equation}
We remark that solution to Equation (\ref{eq:alpha-U-beta}) is
unique because
$\delta^{-1}\{(D+\delta)\alpha(U_\beta)+i/\hbar\Delta_{\alpha\beta}\circ
\alpha(U_\beta)-i/\hbar \alpha(U_\beta)\circ \Delta_\alpha\}$
raises the total degree of $\alpha(U_\beta)$ by 1. Taking $\delta$
on both sides of the above equation, we obtain the first identity
of this Proposition.

(2) Setting $\beta=\alpha^{-1}$ in Equation
(\ref{eq:alpha-U-beta}), we have that
$\alpha(U_{\alpha^{-1}})=1+\delta^{-1}\{(D+\delta)\alpha(U_{\alpha^{-1}})
+i/\hbar(\Delta_{id}\circ \alpha(U_{\alpha^{-1}})-i/\hbar
\alpha(U_{\alpha^{-1}})\circ
\Delta_\alpha\}=1+\delta^{-1}\{(D+\delta)\alpha(U_{\alpha^{-1}})-i/\hbar
\alpha(U_{\alpha^{-1}})\circ \Delta_\alpha\}$, which is same as
the defining equation for $U_\alpha^{-1}$. By the uniqueness of
the solution to the above equation, we have that
$\alpha(U_{\alpha^{-1}})=U_{\alpha}^{-1}$.
\end{proof}

Fedosov \cite{fe} defined quantization of a symplectic
diffeomorphism $\gamma$ on $(M,\omega)$ as
\[
\hat{\gamma}(a)=Ad_{U^{-1}_\gamma}(\gamma(a))=U_\gamma ^{-1}\circ
\gamma(a)\circ U_{\gamma},
\]
which is an algebra endomorphism of the quantum algebra $\calw_D$.

The ``defect" of this quantization is that the associativity of
composition fails, i.e.
\[
\hat{\alpha}\hat{\beta}\ne \widehat{\alpha\beta}.
\]
Instead, Fedosov proves the following property of the associator
$v_{\alpha, \beta}:=U^{-1}_{\alpha}\circ \alpha(U^{-1}_\beta)
\circ \alpha\beta(U^{-1}_{(\alpha\beta)^{-1}})$.

\begin{proposition}
\label{prop:associator} The associator $v_{\alpha,\beta}$ is a
flat section of $\calw$ and
\[
\hat{\alpha}\hat{\beta}\widehat{(\alpha\beta)^{-1}}={\rm{Ad}}_{v_{\alpha,
\beta}}.
\]
\end{proposition}
\begin{proof}
Using $DU^{-1}_\alpha=i/\hbar U^{-1}_\alpha \circ \Delta_\alpha$,
we have
\[
\begin{array}{rccl}
Dv_{\alpha, \beta}&=& & D\Big(U^{-1}_{\alpha}\circ
\alpha(U^{-1}_\beta) \circ
\alpha(\beta(U^{-1}_{(\alpha\beta)^{-1}}))\Big)\cr &=& &
D(U^{-1}_{\alpha})\circ \alpha(U^{-1}_\beta) \circ
\alpha(\beta(U^{-1}_{(\alpha\beta)^{-1}}))\cr & &+&
U^{-1}_{\alpha}\circ D(\alpha(U^{-1}_\beta)) \circ
\alpha\beta(U^{-1}_{(\alpha\beta)^{-1}})\cr & &+&
U^{-1}_{\alpha}\circ \alpha(U^{-1}_\beta) \circ
D(\alpha\beta(U^{-1}_{(\alpha\beta)^{-1}}))\cr &=&
&\displaystyle\frac{i}{\hbar}U^{-1}_{\alpha}\circ
\Delta_\alpha\circ \alpha(U^{-1}_\beta) \circ
\alpha\beta(U^{-1}_{(\alpha\beta)^{-1}})\cr & & -&
\displaystyle\frac{i}{\hbar}U^{-1}_{\alpha}\circ
(\Delta_\alpha\circ\alpha(U^{-1}_\beta)- \alpha(U_\beta
^{-1})\circ \Delta_{\alpha\beta}) \circ
\alpha\beta(U^{-1}_{(\alpha\beta)^{-1}})\cr & &-&
\displaystyle\frac{i}{\hbar}U^{-1}_{\alpha}\circ
\alpha(U^{-1}_\beta) \circ \Delta_{\alpha\beta}\circ
\alpha\beta(U^{-1}_{(\alpha\beta)^{-1}})\cr &=& & 0.
\end{array}
\]
We remark that in the above formula we have used the property
$D(\alpha(U^{-1}_\beta))=i/\hbar \alpha(U_\beta^{-1})\circ
\Delta_{\alpha\beta}-i/\hbar\Delta_{\alpha}\circ
\alpha(U_\beta)^{-1}$. Therefore, by \cite{fe}[Lemma 4.2], we
conclude that $v$ is a flat section of $\calw$.

The property of the associator is a straight forward computation.
\end{proof}

We use $f\mapsto \hat{f}$ to represent the isomorphism between
$C^\infty(M)$ and the quantum algebra $\calw_D$. We know that
$\alpha(\hat{g})$ is a flat section of the connection $\alpha(D)$.
Therefore, $\alpha(\hat{g})$ satisfies the following equation
\[
D(\alpha(\hat{g}))=-\frac{i}{\hbar}[\Delta_\alpha,
\alpha(\hat{g})].
\]
Then $U^{-1}_{\alpha}\circ \alpha(\hat{g})\circ U_{\alpha}$
satisfies the following equation
\[
\begin{split}
&D(U^{-1}_{\alpha}\circ \alpha(\hat{g})\circ U_{\alpha})\\
=&D(U^{-1}_{\alpha})\circ \alpha(\hat{g})\circ U_{\alpha}
+U^{-1}_{\alpha}\circ D(\alpha(\hat{g}))\circ U_{\alpha}
+U^{-1}_{\alpha}\circ \alpha(\hat{g})\circ D(U_{\alpha})\\
=&\frac{i}{\hbar}U_\alpha ^{-1}\circ \Delta_\alpha \circ
\alpha(\hat{g})\circ U_\alpha-\frac{i}{\hbar}U^{-1}_\alpha
\circ\big(\Delta_\alpha\circ \alpha(\hat{g})-\alpha(\hat{g})\circ
\Delta_\alpha\big)\circ U_\alpha-
\frac{i}{\hbar}U_\alpha^{-1}\circ \alpha(\hat{g})\circ\Delta_\alpha\circ U_\alpha\\
=&0
\end{split}
\]

In the following, we apply the above idea to quantize the groupoid
algebra of a pseudogroup $\Gamma$ on a symplectic manifold $M$. We
define the following product on $C_c^\infty(M)\rtimes
\Gamma[[\hbar]]$,
\[
\begin{split}
f \alpha \star g\beta:=&\Big(\hat{f}\circ U^{-1}_\alpha \circ
\alpha(\hat{g})\circ U_\alpha \circ v_{\alpha,
\beta}\Big)\Big|_{y=0}
\alpha\beta\\
=&\Big(\hat{f}\circ U^{-1}_\alpha \circ \alpha(\hat{g})\circ
\alpha(U^{-1}_\beta)\circ
\alpha\beta(U_{(\alpha\beta)^{-1}}^{-1})\Big)\Big|_{y=0}\alpha\beta,
\end{split}
\]
where $y$'s are coordinates along the fiber direction of $T^*M$.
We remark that because $\hat{f}$, $U^{-1}_\alpha \circ
\alpha(\hat{g})\circ U_\alpha$, and $v_{\alpha, \beta}$ are all
flat with respect to the connection $D$, the product $\hat{f}\circ
U^{-1}_\alpha \circ \alpha(\hat{g})\circ \alpha(U^{-1}_\beta)\circ
\alpha\beta(U_{(\alpha\beta)^{-1}}^{-1})$ is also flat with
respect to the connection $D$.

We check the associativity of $\star$ on $C^\infty_c(M)\rtimes
\Gamma[[\hbar]]$. We compute
\[
\begin{split}
&(f\alpha \star g\beta)\star h\gamma\\
=&\Big(\hat{f}\circ U^{-1}_\alpha \circ \alpha(\hat{g})\circ
\alpha(U^{-1}_\beta)\circ
\alpha\beta(U_{(\alpha\beta)^{-1}}^{-1})\Big)|_{y=0}
\alpha\beta\star h\gamma\\
=&\Big(\hat{f}\circ U^{-1}_\alpha \circ \alpha(\hat{g})\circ
\alpha(U^{-1}_\beta)\circ
\alpha\beta(U_{(\alpha\beta)^{-1}}^{-1})\circ
U_{\alpha\beta}^{-1}\circ \alpha\beta(\hat{h})\circ\\
\circ &\alpha\beta(U^{-1}_\gamma)\circ
\alpha\beta\gamma(U^{-1} _{(\alpha\beta\gamma)^{-1}})\Big)|_{y=0}\alpha\beta\gamma\\
=&\Big(\hat{f}\circ U^{-1}_\alpha \circ \alpha(\hat{g})\circ
\alpha(U^{-1}_\beta)\circ \alpha\beta(\hat{h})\circ
\alpha\beta(U^{-1}_\gamma)\circ \alpha\beta\gamma(U^{-1}
_{(\alpha\beta\gamma)^{-1}})\Big)|_{y=0}\alpha\beta\gamma,
\end{split}
\]
where we have used
$\alpha\beta\left(U^{-1}_{(\alpha\beta)^{-1}}\right)=U_{\alpha\beta}$,
and
\[
\begin{split}
&f\alpha\star(g\beta\star h\gamma)\\
=&f\alpha\star\Big(\hat{g}\circ U^{-1}_\beta \circ
\beta(\hat{h})\circ \beta(U_{\gamma} ^{-1})\circ
\beta\gamma(U^{-1}_{(\beta\gamma)^{-1}})
\Big)|_{y=0}\beta\gamma\\
=&\Big(\hat{f}\circ U^{-1}_\alpha\circ \alpha(\hat{g}\circ
U^{-1}_\beta \circ \beta(\hat{h})\circ \beta(U^{-1}_{\gamma})\circ
\beta\gamma(U^{-1}_{(\beta\gamma)^{-1}}))\circ
\alpha(U^{-1}_{\beta\gamma})\circ
\alpha\beta\gamma(U^{-1}_{(\alpha\beta\gamma)^{-1}})\Big)|_{y=0}\alpha\beta\gamma\\
=&\Big(\hat{f}\circ U_\alpha \circ \alpha(\hat{g})\circ
\alpha(U^{-1}_\beta)\circ\alpha\beta(\hat{h})\circ\alpha\beta(U_\gamma^{-1})\circ
\alpha\beta\gamma(U^{-1}_{(\alpha\beta\gamma)^{-1}})\Big)|_{y=0}\alpha\beta\gamma,
\end{split}
\]
where we have used
$\beta\gamma\left(U^{-1}_{(\beta\gamma)^{-1}}\right)=U_{\beta\gamma}$.

We conclude that $\star$ defines an associative product on the
algebra $C^\infty_c(M)\rtimes \Gamma[[\hbar]]$.
\section{A universal deformation formula}
In this section, we apply the construction described in the
previous section to construct a universal deformation formula of
Connes-Moscovici's Hopf algebra $\calh_1$. We start with briefly
recalling the definition of $\calh_1$.

We consider the defining representation of $\calh_1$ on
$C^\infty_c(\reals\times \reals^+)\rtimes \Gamma$. Let $(x,y)$
with $y>0$ be coordinates on $\reals\times \reals^+$. Define
$\gamma:\reals\times \reals^+\to \reals\times \reals^+$ by
\[
\gamma(x,y)=\left(\gamma(x), \frac{y}{\gamma'(x)}\right).
\]
We remark that the above expression of $\gamma$ action does not
agree with the formulas in \cite{CM03-1}, but the two actions are
isomorphic under the transformation $y\mapsto 1/y$.

Consider $X=1/y\partial_x$, and $Y=-y\partial_y$ acting on
$C_c^\infty(\reals\times \reals^+)\rtimes \Gamma$ as follows
\[
X(f\alpha)=\frac{1}{y}f_x\alpha,\ \ Y(f\alpha)=-yf_y\alpha.
\]
We have the following observation
\[
\alpha(Y(\alpha^{-1}(f)))=\alpha\left(Y(f(\alpha(x),
\frac{y}{\alpha'(x)}))\right)=\alpha\left(-\frac{y}{\alpha'(x)}f_y(\alpha(x),
\frac{y}{\alpha'(x)})\right)=-yf_y(x,y)=Y(f),
\]
and
\begin{eqnarray*}
\alpha(X(\alpha^{-1}(f)))&=\alpha\left(X(f(\alpha(x),
\frac{y}{\alpha'(x)}))\right)=\alpha\left(\frac{1}{y}\alpha'(x)f_x(\alpha(x),
\frac{y}{\alpha'(x)})
-\frac{\alpha''(x)}{(\alpha'(x))^2}f_y(\alpha(x),
\frac{y}{\alpha'(x)})\right)\\
&=\frac{1}{y}f_x+\frac{{\alpha^{-1}}''}{{\alpha^{-1}}'}f_y
=Xf-\frac{\log({\alpha^{-1}}')'}{y}Yf=(X-\delta_1(\alpha)Y)f
\end{eqnarray*}
where $\delta_1(f\alpha)=\log({\alpha^{-1}}')'/yf$.

We compute
$\delta_2(\alpha)=X(\delta_1(\alpha))=\frac{1}{y}\partial_x(\log({\alpha^{-1}}')'/y)
=\frac{1}{y^2}\frac{{\alpha^{-1}}'''{\alpha^{-1}}'
-({\alpha^{-1}}'')^2}{({\alpha^{-1}}')^2}$, and define
$\delta_n(\alpha)=X(\delta_{n-1}(\alpha))$ by induction for $n\geq
2$.

On $\reals\times \reals^+$, we consider the Poisson structure
$\partial_x\wedge \partial_y$, which can be expressed by
$-X\otimes Y+Y\otimes X$. Our main goal is to use the method
reviewed in the previous section to construct a star product on
$C^\infty_c(M)\rtimes \Gamma[[\hbar]]$. We prove that this star
product as a bilinear operator actually can be expressed by an
element $R$ of $\calh_1\otimes \calh_1[[\hbar]]$. The
associativity of the star product is equivalent to the property
that $R$ is a universal deformation formula. We start with fixing
a symplectic connection $\nabla$ on the cotangent bundle of
$\reals\times \reals^+$, which was introduced in \cite{bty}[Sec.
3]
\[
\nabla_{\partial_x}\partial_x=0,\nabla_{\partial_x}\partial_y=\frac{1}{2y}\partial_x,
\nabla_{\partial_y}\partial_x=\frac{1}{2y}\partial_x,
\nabla_{\partial_y}\partial_y=-\frac{1}{2y}\partial_y.
\]
Using $X$ and $Y$, we can express the above connection by
\[
\nabla_X X=0,\nabla_X Y=-\frac{1}{2}X, \nabla_Y X=\frac{1}{2}X,
\nabla_YY=-\frac{1}{2}Y.
\]

We compute $\alpha(\nabla)$ by $\alpha \nabla \alpha^{-1}$.
\begin{eqnarray*}
&\alpha(\nabla)_{\partial_x}\partial_x = -\frac{ {\alpha^{-1}}'''
{\alpha^{-1}}'-\frac{3}{2}({\alpha^{-1}}'')^2 }{( {\alpha^{-1}
}')^2}
y\partial_y,&\alpha(\nabla)_{\partial_x}\partial_y=\frac{1}{2y}\partial_x,\\
&\alpha(\nabla)_{\partial_y}\partial_x=\frac{1}{2y}\partial_x,
&\alpha(\nabla)_{\partial_y}\partial_y=-\frac{1}{2y}\partial_y,
\end{eqnarray*}
and
\[
\alpha(\nabla)_X X=\delta_2'(\alpha)Y,\alpha(\nabla)_X
Y=-\frac{1}{2}X, \alpha(\nabla)_Y X=\frac{1}{2}X,
\alpha(\nabla)_YY=-\frac{1}{2}Y,
\]
where $\delta_2'=\delta_2-\frac{1}{2}\delta_1^2$. We observe that
both $\nabla$ and $\alpha(\nabla)$ are flat and torsion free.

We consider the lifting of $\nabla$ and $\alpha(\nabla)$ onto the
Weyl algebra bundle. Use $u,v$ to denote the generators along the
fiber direction of the Weyl algebra bundle $\calw$. We have for
any section $a$ of $\calw$,
\[
\begin{split}
Da&=da-dx \frac{\partial a}{\partial u}-dy\frac{\partial
a}{\partial
v}+\frac{i}{\hbar}[ \frac{1}{2y}v^2 dx+\frac{1}{2y}2uvdy,a],\\
\alpha(D)a&=da-dx \frac{\partial a}{\partial u}-dy\frac{\partial
a}{\partial v}+\frac{i}{\hbar}[
(y^3\delta_2'(\alpha)u^2+\frac{1}{2y}v^2) dx+\frac{1}{2y}2uvdy,a].
\end{split}
\]
Therefore, using the notation of the previous section, we have
that $\Delta_\alpha= y^3 \delta_2'(\alpha)u^2 dx$.

In the following, we solve the expression for $\hat{f},
\alpha(\hat{g}), U^{-1}_\alpha, \alpha(U^{-1}_\beta),
\alpha\beta(U^{-1}_{(\alpha\beta)^{-1}})$.

\subsection{$\hat{f}$}

The section $\hat{f}$ of $\calw$ is a unique solution of
\[
D\hat{f}=0,\ \ \ \hat{f}|_{u=v=0}=f.
\]
Write $\hat{f}=\sum_{mn}f_{mn}u^mv^n$. Then the above equation is
expressed as
\begin{eqnarray*}
&&\sum_{mn}\Big( dx \partial_x f_{mn}u^mv^n+dy \partial_y
f_{mn}u^mv^n-dx f_{mn}mu^{m-1}v^n-dy f_{mn} u^m nv^{n-1} \\
&&+dx \frac{1}{2y}\frac{-1}{2}2vf_{mn}mu^{m-1}v^n+dy
\frac{1}{2y}\frac{1}{2}(2vf_{mn}u^mnv^{n-1}-2uf_{mn}mu^{m-1}v^n)\Big)\\
&&=\sum_{mn}dx\Big( \partial_x
f_{mn}-(m+1)f_{m+1n}-\frac{1}{2y}(m+1)f_{m+1n-1}\Big)u^mv^n\\
&&+\sum_{mn}dy\Big(\partial_yf_{mn}-(n+1)f_{mn+1}+\frac{1}{2y}(n-m)f_{mn}
\Big)u^mv^n.
\end{eqnarray*}
Therefore, we have that $\hat{f}$ is the unique solution to the
following family of equations
 \begin{equation}\label{eq:f}
\begin{split}
&\partial_x f_{mn}-(m+1)f_{m+1n}-\frac{1}{2y}(m+1)f_{m+1n-1}=0\\
&\partial_yf_{mn}-(n+1)f_{mn+1}+\frac{1}{2y}(n-m)f_{mn}=0
\end{split}
\end{equation}
with $f_{00}=f$.

Solving Equation (\ref{eq:f}), we have that
\begin{eqnarray*}
f_{mn}&=&\frac{1}{m!n!}\left(\partial_y-\frac{m+1-n}{2y}\right)\cdots\left(\partial_y-\frac{m}{2y}\right)\partial_x^mf\\
&=&\frac{(-1)^{n}}{m!}y^{m-n}X^m\left(Y+\frac{m+n-1}{2}\right)\cdots\left(Y+\frac{m}{2}\right)(f).
\end{eqnarray*}
\subsection{$\alpha(\hat{g})$}

We know that $\alpha(\hat{g})$ is the unique solution to the
following equation
\[
D(\alpha(\hat{g}))=-\frac{i}{\hbar}[\Delta_\alpha,
\alpha(\hat{g})],
\]
with $\alpha(\hat{g})|_{u=v=0}=\alpha(g)$. We remind that
$\Delta_\alpha=y^3\delta_2'(\alpha)u^2dx$.

Similar to $\hat{f}$, $\alpha(\hat{g})$ satisfies the following
equation
\[
\begin{split}
0=&\sum_{mn}dx\Big( \partial_x
\alpha(g)_{mn}-(m+1)\alpha(g)_{m+1n}\\
&-\frac{1}{2y}(m+1)\alpha(g)_{m+1n-1}
+y^3\delta_2'(n+1)\alpha(g)_{m-1n+1} \Big)u^mv^n\\
&+\sum_{mn}dy\Big(\partial_y\alpha(g)_{mn}-(n+1)\alpha(g)_{mn+1}+\frac{1}{2y}(n-m)\alpha(g)_{mn}
\Big)u^mv^n
\end{split}
\]
with $\alpha(g)_{00}=\alpha(g)$. Therefore, $\alpha(\hat{g})$ is
the unique solution to the following family of equations
\begin{equation}\label{eq:g}
\begin{split}
0&=\partial_x \alpha(g)_{mn}-(m+1)\alpha(g)_{m+1n}
-\frac{1}{2y}(m+1)\alpha(g)_{m+1n-1}
+y^3\delta_2'(n+1)\alpha(g)_{m-1n+1}\\
0&=\partial_y\alpha(g)_{mn}-(n+1)\alpha(g)_{mn+1}+\frac{1}{2y}(n-m)\alpha(g)_{mn}
\end{split}
\end{equation}
with $\alpha(g)_{00}=\alpha(g)$.

By the second equation of Equations (\ref{eq:g}), we have
\[
\alpha(g)_{m,n}=\frac{1}{n}\left(\partial_y+\frac{n-1-m}{2y}\right)\alpha(g)_{mn-1}
=\cdots=\frac{1}{n!}\left(\partial_y+\frac{n-1-m}{2y}\right)\cdots\left(\partial_y+\frac{-m}{2y}\right)\alpha(g)_{m,0}.
\]

Setting $n=0$ in the first equation of Equations (\ref{eq:g}), we
obtain
\[
\begin{split}
\alpha(g)_{m+1,0}&=\frac{1}{m+1}(\partial_x
\alpha(g)_{m,0}+y^3\delta_2'\alpha(g)_{m-1,1})\\
&=\frac{1}{m+1}\left(\partial_x
\alpha(g)_{m,0}+y^3\delta_2'\left(\partial_y+\frac{1-m}{2y}\right)\alpha(g)_{m-1,0}\right).
\end{split}
\]

By induction, we can solve the above equation as
\[
\alpha(g)_{mn}=\frac{(-1)^ny^{m-n}}{m!n!}A_m\left(Y+\frac{n+m-1}{2}\right)\cdots\left(Y+\frac{m}{2}\right)(\alpha(g)),
\]
where $A_m\in \calh_1$ is defined inductively by
\[
A_{m+1}=XA_m-m\delta_2'\left(Y-\frac{m-1}{2}\right)A_{m-1},\ \
A_0=1.
\]

\subsection{$U_\alpha^{-1}$}

We compute $U^{-1}_\alpha$ using the equation
\[
DU^{-1}_{\alpha}=U^{-1}_\alpha\circ \frac{i}{\hbar}\Delta_\alpha
\]
with $\Delta_\alpha=\delta_2'(\alpha)y^3u^2 dx$.

Write $U^{-1}_{\alpha}=\sum_{mn}u^\alpha_{mn}u^mv^n$, where
$u^\alpha_{mn}$ takes values in $C^\infty_c(M)[\hbar^{-1},
\hbar]]$. And $u^\alpha_{mn}$ satisfies the following family of
equations
\begin{equation}\label{eq:u-ind}
\begin{split}
0&=\partial_x u^\alpha_{mn}-(m+1)u^\alpha_{m+1n}
-\frac{1}{2y}(m+1)u^\alpha_{m+1n-1}
-\frac{i}{\hbar}y^3\delta_2'u^\alpha_{m-2n}\\
&+y^3\delta_2'(n+1)u^\alpha_{m-1n+1}
+\frac{i\hbar}{4}y^3\delta_2'(n+2)(n+1)u^\alpha_{mn+2},\\
0&=\partial_yu^\alpha_{mn}-(n+1)u^\alpha_{mn+1}+\frac{1}{2y}(n-m)u^\alpha_{mn}
\end{split}
\end{equation}
with $u_{0,0}=1$.

The second equation of Equations (\ref{eq:u-ind}) implies that
\[
u^\alpha_{mn}=\frac{1}{n}(\partial_y+\frac{n-1-m}{2y})u^\alpha_{mn-1}
=\cdots=\frac{1}{n!}(\partial_y+\frac{n-1-m}{2y})\cdots
(\partial_y+\frac{-m}{2y})u^\alpha_{m,0}.
\]

We use the $n=0$ version of the first equation of Equation
(\ref{eq:u-ind}) to solve $u_{m0}$.
\[
\begin{split}
u^\alpha_{m+10}&=\frac{1}{m+1}\Big(\partial_xu^\alpha_{m0}
-\frac{i}{\hbar}y^3\delta_2'u^\alpha_{m-20}+y^3\delta_2'(\partial_y-\frac{m-1}{2y})
u^\alpha_{m-1}\\
&\hspace{1.5cm}
+\frac{i\hbar}{4}y^3\delta_2'(\partial_y-\frac{m-1}{2y})(\partial_y-\frac{m}{2y})u^\alpha_{m0}\Big)
\end{split}
\]

By induction, we have the following expression of $u$,
\[
u^\alpha_{mn}=\frac{(-1)^ny^{m-n}}{m!n!}(Y+\frac{n-m-1}{2})\cdots(Y-\frac{m}{2})A_{m}1,
\]
where $A_m$ is defined by
\[
A_{m+1}=\Big(X+\frac{i\hbar}{4}\delta_2'(Y-\frac{m-1}{2})(Y-\frac{m}{2})\Big)A_m
-\delta_2'(Y-\frac{m-1}{2})A_{m-1}-\frac{i}{\hbar}\delta_2'A_{m-2},
\]
with $A_0=1$.
\begin{remark}
We need to prove that the above obtained solution
$\tilde{U}_\alpha^{-1}=\sum u^\alpha_{mn}u^mv^n$ is the unique
solution to the defining Equation (\ref{eq:u-inverse}) of
$U_\alpha^{-1}$, which implies that the above
$\tilde{U}_\alpha^{-1}=U_\alpha^{-1}$.

Using the above expression for $u^\alpha_{mn}$, by induction we
can prove that the negative power of $\hbar$ contained in
$u^\alpha_{mn}$ is less than or equal to $[m/3]$ ($[\mu]$ means
the Gauss integer function.). This shows that once $m+n>0$, the
smallest degree term contained in $u^\alpha_{mn}u^m v^n$ has
degree greater than or equal to 1. Therefore the degree 0 term of
the solution $\tilde{U}_\alpha^{-1}=\sum u^\alpha_{mn}u^mv^n$ is
equal to 1. Accordingly, we compute using that
$D\tilde{U}_\alpha^{-1}-\tilde{U}_\alpha^{-1}\circ i/\hbar
\Delta_\alpha=0$,
\begin{eqnarray*}
\tilde{U}^{-1}_\alpha &=& \delta
\delta^{-1}\tilde{U}^{-1}_\alpha+\delta^{-1}\delta\tilde{U}^{-1}_\alpha+1=1+
\delta^{-1}\delta\tilde{U}^{-1}_\alpha\\
&=&1+\delta^{-1}(\delta\tilde{U}^{-1}_\alpha+D\tilde{U}_\alpha^{-1}-\tilde{U}_\alpha^{-1}\circ
i/\hbar \Delta_\alpha)\\
&=&1+\delta^{-1}\{(D+\delta)\tilde{U}^{-1}_\alpha-\tilde{U}^{-1}_\alpha\circ
i/\hbar\Delta_\alpha\}.
\end{eqnarray*}

This remark applies also to the following solutions
$\alpha(U_\beta^{-1})$ and
$\alpha\beta(U^{-1}_{(\alpha\beta)^{-1}})$.
\end{remark}
\subsection{$\alpha(U_\beta^{-1})$ and
$\alpha\beta(U^{-1}_{(\alpha\beta)^{-1}})$}

By Proposition \ref{prop:u}, we know that $\alpha(U_\beta)$
satisfies the equation
\[
D(\alpha(U_\beta))=-\frac{i}{\hbar}\big(\Delta_{\alpha\beta}\circ
\alpha(U_{\beta})-\alpha(U_\beta)\circ \Delta_{\alpha}\big).
\]
Accordingly, $\alpha(U_\beta^{-1})=(\alpha(U_\beta))^{-1}$
satisfies
\[
\begin{split}
D(\alpha(U_\beta^{-1}))&=-(\alpha(U_\beta))^{-1}\circ
D(\alpha(U_\beta))\circ (\alpha(U_\beta))^{-1}\\
&=-(\alpha(U_\beta))^{-1}\circ\big(-\frac{i}{\hbar}\Delta_{\alpha\beta}\circ
\alpha(U_{\beta})+\frac{i}{\hbar} \alpha(U_\beta)\circ
\Delta_{\alpha}\big)\circ (\alpha(U_\beta))^{-1}\\
&=\frac{i}{\hbar}\big(\alpha(U_\beta^{-1})\circ
\Delta_{\alpha\beta}-\Delta_{\alpha}\circ
\alpha(U^{-1}_{\beta})\big)\\
\end{split}
\]

If we write $\alpha(U_{\beta}^{-1})=\sum
u^{\alpha,\beta}_{mn}u^mv^n$, then we have
\begin{equation}\label{eq:u-alpha-beta}
\begin{split}
0&=\partial_xu^{\alpha,\beta}_{mn}-(m+1)u^{\alpha,\beta}_{m+1n}-\frac{1}{2y}(m+1)
u^{\alpha,\beta}_{m+1n-1}-\frac{i}{\hbar}y^3\alpha(\delta_2'(\beta))u^{\alpha,\beta}_{m-2n}\\
&+y^3(n+1)(2\delta_2'(\alpha)+\alpha(\delta_2'(\beta)))u^{\alpha,\beta}_{m-1n+1}
+\frac{i\hbar}{4}y^3\alpha(\delta'_2(\beta))u^{\alpha,\beta}_{mn+2}\\
0&=\partial_yu^{\alpha,\beta}_{mn}-(n+1)u^{\alpha,\beta}_{mn+1}+\frac{1}{2y}(n-m)
u^{\alpha,\beta}_{mn}.
\end{split}
\end{equation}

We can solve Equation (\ref{eq:u-alpha-beta}) of
$u^{\alpha,\beta}_{mn}$ as follows
\[
u^{\alpha,\beta}_{m,n}=\frac{(-1)^ny^{m-n}}{m!n!}(Y+\frac{n-m-1}{2})\cdots(Y-\frac{m}{2})A_{m}1,
\]
where $A_m\in \calh_1$ is defined inductively
\[
\begin{split}
A_{m+1}&=\left(X+\frac{i\hbar}{4}\alpha(\delta_2'(\beta))(Y-\frac{m-1}{2})(Y-\frac{m}{2})\right)A_m\\
&-\big(2\delta_2'(\alpha)
+\alpha(\delta_2'(\beta))\big)(Y-\frac{m-1}{2})A_{m-1}-
\frac{i}{\hbar}\alpha(\delta_2'(\beta))A_{m-2}, \ \ \ A_0=1.
\end{split}
\]

We know that $\alpha\beta(U^{-1}_{(\alpha\beta)^{-1}})$ is equal
to $U_{\alpha\beta}$, which satisfies
\[
DU_{\alpha\beta}=-\frac{i}{\hbar}\Delta_{\alpha\beta}\circ
U_{\alpha\beta}.
\]
We can solve $\alpha\beta(U^{-1}_{(\alpha\beta)^{-1}})$ as
$U_{\alpha\beta}$. Write
$\alpha\beta(U^{-1}_{(\alpha\beta)^{-1}})=U_{\alpha\beta}=\sum
v^{\alpha\beta}_{mn}u^mv^n$. Then
\[
v^{\alpha\beta}_{mn}=\frac{(-1)^ny^{m-n}}{m!n!}(Y+\frac{n-m-1}{2})\cdots(Y-\frac{m}{2})A_{m}1,
\]
where $A_m\in \calh_1$ is defined by
\[
A_{m+1}=\left(X-\frac{i\hbar}{4}\delta_2'(\alpha\beta)(Y-\frac{m-1}{2})(Y-\frac{m}{2})\right)A_m
+\delta_2'(\alpha\beta)(Y-\frac{m-1}{2})A_{m-1}+\frac{i}{\hbar}\delta_2'(\alpha\beta)A_{m-2},
\]
with $A_0=1$.

We notice that terms surviving in the product $\hat{f}\circ
U_\alpha^{-1}\circ \alpha(\hat{g})\circ \alpha(U_\beta^{-1})\circ
\alpha\beta(U^{-1}_{(\alpha\beta)^{-1}})|_{u=v=0}$ are sums of
terms of the following form
\[
C_{m_1,\cdots,m_5;n_1,\cdots,
n_5}=f_{m_1n_1}u^\alpha_{m_2n_2}\alpha(g)_{m_3n_3}u^{\alpha,\beta}_{m_4n_4}v^{\alpha\beta}_{m_5n_5}
u^{m_1}v^{n_1}\circ \cdots \circ u^{m_5}v^{n_5}|_{u=v=0}
\]
with $m_1+\cdots+m_5=n_1+\cdots+n_5$.

\begin{theorem}\label{thm:main}
There exists an element $R\in \calh_1\otimes \calh_1[[\hbar]]$
such that the star product on
$C^\infty_c(\reals\times\reals^+)\rtimes \Gamma$ can be expressed
by $f\alpha\star g\beta=m(R(f\alpha\otimes g\beta))$, where
$m:C_c^\infty(\reals\times \reals^+)\rtimes \Gamma\otimes
C_c^\infty(\reals\times \reals^+)\rtimes \Gamma \to
C_c^\infty(\reals\times \reals^+)\rtimes \Gamma$ is the
multiplication map. Furthermore, $R$ is a universal deformation
formula of $\calh_1$.
\end{theorem}
\begin{proof}
We start by rewriting
$u^\alpha_{mn},u^{\alpha,\beta}_{mn},v^{\alpha\beta}_{mn}$, and
$\alpha(\hat{g})$.
\begin{enumerate}
\item As $X$ and $Y$ vanishes on 1, $u^\alpha_{mn}$ can be
written as $y^{m-n}$ times a sum of terms like powers of $X$ and
$Y$ acting on powers of $\delta_2'(\alpha)$. If we rewrite
$\delta_2'(\alpha)$ as $\delta_2-1/2\delta_1^2$, we can write a
term of powers of $X$ and $Y$ acting on powers of
$\delta_2'(\alpha)$ into a sum of products like
$\delta_{i_1}^{j_1}(\alpha)\cdots\delta_{i_p}^{j_p}(\alpha)$. We
point out that as there are $1/\hbar$ in the induction formula of
$A_m$, $u^\alpha_{mn}$ may contain terms with negative powers of
$\hbar$. However, by induction on the negative power of $\hbar$,
we can easily prove that the negative power of $\hbar$ is no more
than $[m/3]$. Therefore, we can write $u^\alpha_{mn}$ as
$\hbar^{-[\frac{m}{3}]}y^{m-n}\mu_{mn}(\delta_1(\alpha),\delta_2(\alpha),\cdots)$,
where $\mu_{mn}$ is a polynomial of variables $\hbar, \delta_1,
\delta_2, \cdots $ independent of $\alpha$.
\item Analogous to the above analysis, we have that
$u^{\alpha,\beta}_{m,n}$ can be written as $y^{m-n}$ times a sum
of terms like powers of $X$ and $Y$ acting on products of powers
of $\delta'_2(\alpha)$ and $\alpha(\delta_2'(\beta))$. When $X$
and $Y$ act on $\delta_2'(\alpha)$, we can write the resulting
terms into polynomials of $\delta_1(\alpha),\cdots,
\delta_p(\alpha),\cdots$. To compute $X$ and $Y$ action on
$\alpha(\delta_2'(\beta))$, we look at the following properties of
$X$ and $Y$, for any function $f$
\[
\begin{split}
X(\alpha(f))&=\alpha(\alpha^{-1}(X(\alpha(f)))
=\alpha(X(f)-\delta_1(\alpha^{-1})Y(f))\\
&=\alpha(X(f))-\alpha(\delta_1(\alpha^{-1}))\alpha(Y(f))=\alpha(X(f))+\delta_1(\alpha)\alpha(Y(f))\\
Y(\alpha(f))&=\alpha(\alpha^{-1}(Y(\alpha(f)))) =\alpha(Y(f)),
\end{split}
\]
where we have used the commutation relation between $X$, $Y$ and
$\alpha$. This implies that powers of $X$, $Y$ acting on
$\alpha(\delta_2'(\beta))$ gives a sum of terms like
$\sigma(\delta_1(\alpha),\cdots)\alpha(\xi(\delta_1(\beta),\cdots))$
with $\nu,\xi$ polynomials in $\calh_1[\hbar]$ independent of
$\alpha,\beta$. We summarize that $u^{\alpha,\beta}_{mn}$ can be
written as
\[
\hbar^{-[\frac{m}{3}]}y^{m-n}\sum_i\nu^i_{mn}(\delta_1(\alpha),\cdots)\alpha(\xi^i_{mn}(\delta_1(\beta),\cdots))
\]
where $\nu^i,\xi^i$ are polynomials of $\hbar, \delta_1, \delta_2,
\cdots,$ independent of $\alpha,\beta$.
\item Similar to $u^\alpha_{mn}$, $v^{\alpha\beta}_{mn}$ can be
written as a sum of terms like powers of $X,Y$ acting on
$\delta_2'(\alpha\beta)$. For our purpose, we need to rewrite
$\delta_2'(\alpha\beta)$ by a sum of
$\delta_2'(\alpha)+\alpha(\delta_2'(\beta))$. Therefore the
situation is similar to $u^{\alpha,\beta}_{mn}$. We can write
$v^{\alpha\beta}_{mn}$ as
\[
\hbar^{-[\frac{m}{3}]}y^{m-n}\sum
\eta^j_{mn}(\delta_1(\alpha),\cdots)\alpha(\lambda^j_{mn}(\delta_1(\beta),\cdots))
\]
with $\eta^j_{mn},\lambda^j_{mn}$ polynomials independent of
$\alpha,\beta$.
\item From the inductive relations, we see that $\alpha(\hat{g})_{mn}$ can be written as a sum
of terms like product of two parts. One part is powers of $X$ and
$Y$ acting on $\delta_2'(\alpha)$, the other is powers of $X$ and
$Y$ acting on $\alpha(g)$. We can write the part involving
$\delta_2'(\alpha)$ as polynomials of
$\delta_1(\alpha),\delta_2(\alpha),\cdots$, the part with
$\alpha(g)$ like the above $\alpha(\delta_2'(\beta))$ as a sum of
terms like
\[
\varphi(\delta_1(\alpha),\cdots)\alpha(\phi(X,Y)(g)).
\]
Therefore, we can write $\alpha(\hat{g})_{mn}$ as
\[
\sum
y^{m-n}\rho^i_{mn}(\delta_1(\alpha),\cdots)\alpha(\psi^i_{mn}(X,Y)(g)).
\]
\end{enumerate}

Summarizing the above consideration, we can write the term
$C_{m_1,\cdots, m_5;n_1,\cdots,n_5}$ as
\[
\begin{split}
c_{m_1,\cdots,m_5;n_1,\cdots,n_5}\cdot&\hbar^{m_1+m_2-[\frac{m_2}{3}]+m_3+m_4-[\frac{m_4}{3}]
+m_5-[\frac{m_5}{3}]}\sum_{i,j,k}\tau_{m_1n_1}(X,Y)(f)\mu_{m_2n_2}(\delta_1(\alpha),\cdots)\\
&\rho^{i}_{m_3n_3}(\delta_1(\alpha),\cdots)\alpha(\psi^i_{m_3n_3}(X,Y)(g))\\
&v^j_{m_4n_4}(\delta_1(\alpha),\cdots)\alpha(\xi^j_{m_4n_4}(\delta_1(\beta),\cdots))
\eta^{k}_{m_5n_5}(\delta_1(\alpha),\cdots)\alpha(\lambda^k_{m_5n_5}(\delta_1(\beta),\cdots)),
\end{split}
\]
where $c_{m_1,\cdots,m_5;n_1,\cdots,n_5}$ is a constant. And
$\hat{f}\circ U_\alpha^{-1}\circ \alpha(\hat{g})\circ
\alpha(U_\beta^{-1})\circ
\alpha\beta(U^{-1}_{(\alpha\beta)^{-1}})|_{u=v=0}\alpha\beta$ can
be written
\begin{eqnarray*}
&&\sum_{\begin{array}{l}m_1,\cdots, m_5;n_1,\cdots, n_5\\
m_1+\cdots+m_5=n_1+\cdots+n_5
\end{array}}c_{m_1,\cdots,m_5;n_1,\cdots,n_5}
\hbar^{m_1+m_2-[\frac{m_2}{3}]+m_3+m_4-[\frac{m_4}{3}]+m_5-[\frac{m_5}{3}]}\\
&&\sum_{i,j,k}\tau_{m_1n_1}(X,Y)(f)\mu_{m_2n_2}(\delta_1(\alpha),\cdots)
\rho^{i}_{m_3n_3}(\delta_1(\alpha),\cdots)\alpha(\psi^i_{m_3n_3}(X,Y)(g))\\
&&v^j_{m_4n_4}(\delta_1(\alpha),\cdots)\alpha(\xi^j_{m_4n_4}(\delta_1(\beta),\cdots))
\eta^k_{m_5n_5}(\delta_1(\alpha),\cdots)\alpha(\lambda^k_{m_5n_5}(\delta_1(\beta),\cdots))\alpha\beta\\
&=&\sum_{m_1,\cdots, m_5;n_1,\cdots,
n_5}c_{m_1,\cdots,m_5;n_1,\cdots,n_5}\hbar^{m_1+m_2-[\frac{m_2}{3}]+m_3+m_4-[\frac{m_4}{3}]
+m_5-[\frac{m_5}{3}]}\\
&&\sum_{i,j,k}\mu_{m_2n_2}(\delta_1(\alpha),\cdots)\rho^{i}_{m_3n_3}(\delta_1(\alpha),\cdots)
v^j_{m_4n_4}(\delta_1(\alpha),\cdots)\eta^k_{m_5n_5}(\delta_1(\alpha),\cdots)\tau_{m_1n_1}(X,Y)(f)\\
&&\alpha\Big(\xi^j_{m_4n_4}(\delta_1(\beta),\cdots)\lambda^k_{m_5n_5}(\delta_1(\beta),\cdots)
\psi^i_{m_3n_3}(X,Y)(g)\Big)\alpha\beta
\end{eqnarray*}
\begin{eqnarray*}
&=&\sum_{m_1,\cdots, m_5;n_1,\cdots,
n_5}c_{m_1,\cdots,m_5;n_1,\cdots,n_5}\hbar^{m_1+m_2-[\frac{m_2}{3}]+m_3+m_4-[\frac{m_4}{3}]
+m_5-[\frac{m_5}{3}]}\\
&&\sum_{i,j,k}\mu_{m_2n_2}(\delta_1(\alpha),\cdots)\rho^{i}_{m_3n_3}(\delta_1(\alpha),\cdots)
v^j_{m_4n_4}(\delta_1(\alpha),\cdots)\eta^k_{m_5n_5}(\delta_1(\alpha),\cdots)\tau_{m_1n_1}(X,Y)(f)\alpha \\
&&\xi^j_{m_4n_4}(\delta_1(\beta),\cdots)\lambda^k_{m_5n_5}(\delta_1(\beta),\cdots)
\psi^i_{m_3n_3}(X,Y)(g)\beta
\end{eqnarray*}

Define $R_{m_1,\cdots, m_5;n_1,\cdots, n_5}\in
\calh_1[\hbar]\otimes \calh_1[\hbar]$ as
\[
\begin{split}
c_{m_1,\cdots,m_5;n_1,\cdots,n_5}\sum_{i,j,k}&\mu_{m_2n_2}(\delta_1,\cdots)
\rho^i_{m_3n_3}(\delta_1,\cdots)\nu^j_{m_4n_4}(\delta_1,\cdots)
\eta^k_{m_5n_5}(\delta_1,\cdots)\tau_{m_1n_1}(X,Y)\\
\otimes&\xi^j_{m_4n_5}(\delta_1,\cdots)\lambda^k_{m_5n_5}(\delta_1,\cdots)\psi^i_{m_3n_3}(X,Y)
\end{split}
\]
Furthermore, we define
\[
R=\sum_{m_1+\cdots+m_5=n_1+\cdots+n_5}
\hbar^{m_1+m_2-[\frac{m_2}{3}]+m_3+m_4-[\frac{m_4}{3}]+m_5-[\frac{m_5}{3}]}R_{m_1,\cdots,m_5;n_1,\cdots,n_5}.
\]

We conclude that $R\in \calh_1[[\hbar]]\otimes \calh_1[[\hbar]]$
satisfies
\[
f\alpha\star g\beta=m(R(f\alpha\otimes g\beta)).
\]

To check that $R$ is a universal deformation formula, we need to
check that
\begin{enumerate}
\item $(\Delta\otimes 1)R(R\otimes1)=(1\otimes \Delta)R(1\otimes
R)$. This is from the fact that $\star$ is associative and the
$\calh_1$ action on the collection of $C^\infty_c(\reals\times
\reals^+)\rtimes \Gamma$ over all pseudogroups $\Gamma$ is fully
injective.
\item $(\epsilon\otimes 1)R=1\otimes 1=(1\otimes \epsilon)R$. This
is equivalent to show that 1 is a unit respect to the $\star$
product on $C_c^\infty(\reals\times\reals^+)\rtimes \Gamma$.

When $f$ is 1 and $\alpha$ is identity, we have that $\hat{f}=1$,
$U_{\rm {id}}^{-1}=1$, and
$U^{-1}_\beta\circ\beta(U^{-1}_{\beta^{-1}})=1$. This implies that
$1\star g\beta=g\beta$.

When $g$ is 1 and $\beta$ is identity, we have that
$\alpha(\hat{1})=1$, $U^{-1}_\beta=1$ and $U_\alpha^{-1}\circ
\alpha(U^{-1}_{\alpha^{-1}})=1$. Therefore, $f\alpha\star
1=f\alpha$.
\end{enumerate}
\end{proof}

From the computation in the next section, we know that $R\in
1\otimes 1+\hbar R'$, where $R'$ is an element in
$\calh_1[[\hbar]]\otimes \calh_1[[\hbar]]$. Therefore, $R$ is
invertible with $R^{-1}=1+\sum_i(-1)^i(\hbar R')^i$. By the
property of universal deformation formula, we can introduce a new
Hopf algebra structure on $\calh_1[[\hbar]]$ by twisting the
coproduct $\Delta$ by
\[
\tilde{\Delta}(a)=R^{-1}\Delta(a)R,
\]
and the antipode by
\[
\tilde{S}(a)=v^{-1}S(a)v,
\]
with $v=m(S\otimes 1)(R)$.
\section{Formulas of lower order terms}
We must say that the formula for $R$ constructed in the previous
section Theorem \ref{thm:main} could be very complicated. We do
not know an easy way to write it down explicitly. We will provide
in this section the computation of $R$ up to the second order of
$\hbar$.

It is not difficult to check that if
$m_1+m_2-[m_2/3]+m_3+m_4-[m_4/3]+m_5-[m_5/3]=0$ for nonnegative
integers $m_i$, $i=1,\cdots,5$, then $m_1=\cdots=m_5=0$.
Therefore, the $\hbar^0$ component of the $\star$ product
$f\alpha\star g\beta$ is equal to $f\alpha(g)\alpha\beta$. This
implies that $R_{0,\cdots,0;0,\cdots,0}=1$ and $R=1\otimes
1+O(\hbar)$.

We consider $R_{m_1,\cdots,m_5;n_1,\cdots,n_5}$ with
$m_1+m_2-[m_2/3]+m_3+m_4-[m_4/3]+m_5-[m_5/3]=1$. It is not
difficult to find that this implies one of $m_i$, $i=1,\cdots, 5$
takes value 1, and all others vanish. We also check that
$u^\alpha_{11}=u^\alpha_{01}=u^{\alpha}_{10}=u^{\alpha,\beta}_{11}
=u^{\alpha,\beta}_{10}=u^{\alpha,\beta}_{01}=v^{\alpha\beta}_{11}
=v^{\alpha\beta}_{10}=v^{\alpha\beta}_{01}=0$. Therefore,
$R_{1,0,\cdots,0;0,0,1,0,0}$ and $R_{0,0,1,0,0;1,0,\cdots,0}$ are
the only nonzero terms among all
$R_{m_1,\cdots,m_5;n_1,\cdots,n_5}$ with
$$m_1+m_2-[m_2/3]+m_3+m_4-[m_4/3]+m_5-[m_5/3]=1.$$
We compute $R_{1,0,\cdots,0;0,0,1,0,0}=\frac{i\hbar}{2}X\otimes
Y$, and
$R_{0,0,1,0,0;1,0,\cdots,0}=-\frac{i\hbar}{2}(\delta_1Y\otimes
Y+Y\otimes X)$. Therefore, the $\hbar$ component of $R$ is
\[
\frac{-i\hbar}{2}(-X\otimes Y+\delta_1Y\otimes Y+Y\otimes
X)=\frac{-i\hbar}{2}(S(X)\otimes Y+Y\otimes X).
\]

We consider $R_{m_1,\cdots,m_5;n_1,\cdots, n_5}$ with
$m_1+m_2-[m_2/3]+m_3+m_4-[m_4/3]+m_5-[m_5/3]=2$. There are three
classes of possibilities. i) one of $m_2,m_4, m_5$ is equal to 3;
ii)one of $m_i$ ($i=1,\cdots,5$) is equal to $2$, iii) two of
$m_i$ ($i=1,\cdots,5$) are both equal to 1. We notice that
$u^\alpha_{1\cdot}=u^\alpha_{2\cdot}=u^{\alpha,\beta}_{1\cdot}=u^{\alpha,\beta}_{2\cdot}
=v^{\alpha\beta}_{1\cdot}=v^{\alpha\beta}_{2\cdot}=0$. This
implies that terms contributing in
$m_1+m_2-[m_2/3]+m_3+m_4-[m_4/3]+m_5-[m_5/3]=2$ are from following
three groups.
\begin{enumerate}
\item $R_{0,3,0,0,0;1,0,2,0,0}, R_{0,3,0,0,0;2,0,1,0,0},
R_{0,0,0,3,0;1,0,2,0,0}$,\\
$R_{0,0,0,3,0;2,0,1,0,0}, R_{0,0,0,0,3;1,0,2,0,0}$,
$R_{0,0,0,0,3;2,0,1,0,0}$;
\item $R_{0,3,0,0,0;3,0,\cdots,0}$, $R_{0,3,0,0,0;0,0,3,0,0}$,
$R_{0,0,0,3,0;3,0,\cdots,0}$,\\
$R_{0,0,0,3,0;0,0,3,0,0}$, $R_{0,0,0,0,3;3,0,\cdots,0}$,
$R_{0,0,0,0,3;0,0,3,0,0}$;
\item $R_{2,0,\cdots,0;0,0,2,0,0}, R_{0,0,2,0,0;2,0,0,0,0}$, $R_{1,0,1,0,0;1,0,1,0,0}$.
\end{enumerate}

We compute the above terms separately.
\begin{itemize}
\item $R_{0,3,0,0,0;1,0,2,0,0}=(\frac{-i\hbar}{2})^{3}\frac{1}{2}
\frac{-i}{\hbar}\delta_2'Y\otimes(Y+\frac{1}{2})Y=-(\frac{-i\hbar}{2})^2\frac{1}{4}\delta_2'Y\otimes(Y+\frac{1}{2})Y
$.

\item $R_{0,3,0,0,0;2,0,1,0,0}=(\frac{-i\hbar}{2})^3\frac{1}{2!}\frac{i}{\hbar}
\delta_2'(Y+\frac{1}{2})Y\otimes Y
=(\frac{-i\hbar}{2})^2\frac{1}{4}\delta_2'(Y+\frac{1}{2})Y\otimes
Y$.

\item $R_{0,0,0,3,0;1,0,2,0,0}=-(\frac{-i\hbar}{2})^3\frac{1}{2}\frac{i}{\hbar}Y
\otimes
\delta_2'(Y+\frac{1}{2})Y=-(\frac{-i\hbar}{2})^2\frac{1}{4}Y
\otimes \delta_2'(Y+\frac{1}{2})Y$.

\item $R_{0,0,0,3,0;2,0,1,0,0}=
-(\frac{-i\hbar}{2})^3\frac{1}{2}\frac{i}{\hbar}(Y+\frac{1}{2})Y\otimes
\delta_2'
Y=-(\frac{-i\hbar}{2})^2\frac{1}{4}(Y+\frac{1}{2})Y\otimes\delta_2'Y$.

\item $R_{0,0,0,0,3;1,0,2,0,0}=(\frac{-i\hbar}{2})^3\frac{1}{2}\frac{i}{\hbar}\big[\delta_2'Y\otimes
(Y+\frac{1}{2})Y+Y\otimes \delta_2'
(Y+\frac{1}{2})Y\big]=(\frac{-i\hbar}{2})^2\frac{1}{4}\big[\delta_2'Y\otimes
(Y+\frac{1}{2})Y+Y\otimes \delta_2' (Y+\frac{1}{2})Y\big]$.

\item $R_{0,0,0,0,3;2,0,1,0,0}=(\frac{-i\hbar}{2})^3\frac{1}{2}\frac{i}{\hbar}
\big[\delta_2'(Y+\frac{1}{2})Y\otimes Y+(Y+\frac{1}{2})Y\otimes
\delta_2'Y\big]=(\frac{-i\hbar}{2})^2\frac{1}{4}\big[\delta_2'(Y+\frac{1}{2})Y\otimes
Y+(Y+\frac{1}{2})Y\otimes \delta_2'Y\big]$.

\item $R_{0,3,0,0,0;3,0,0,0,0}=-(\frac{-i\hbar}{2})^3\frac{i}{6\hbar}
\delta_2'(Y+1)(Y+\frac{1}{2})Y\otimes1=-(\frac{-i\hbar}{2})^2\frac{1}{12}\delta_2'(Y+1)
(Y+\frac{1}{2})Y\otimes 1$.

\item $R_{0,3,0,0,0;0,0,3,0,0}=
(\frac{-i\hbar}{2})^3\frac{i}{6\hbar}\delta_2'\otimes(Y+1)(Y+\frac{1}{2})Y
=(\frac{-i\hbar}{2})^2\frac{1}{12}\delta_2'\otimes
(Y+1)(Y+\frac{1}{2})Y$.

\item $R_{0,0,0,3,0;3,0,0,0,0}=-(\frac{-i\hbar}{2})^3\frac{i}{6\hbar}
(Y+1)(Y+\frac{1}{2})Y\otimes\delta_2'=-(\frac{-i\hbar}{2})^2\frac{1}{12}
(Y+1)(Y+\frac{1}{2})Y\otimes\delta_2'$.

\item $R_{0,0,0,3,0;0,0,3,0,0}=-(\frac{-i\hbar}{2})^3\frac{i}{6\hbar}
1\otimes\delta_2'(Y+1)(Y+\frac{1}{2})Y=-(\frac{-i\hbar}{2})^2\frac{1}{12}
1\otimes\delta_2'(Y+1)(Y+\frac{1}{2})Y$.

\item $R_{0,0,0,0,3;3,0,0,0,0}=(\frac{-i\hbar}{2})^3\frac{i}{6\hbar}
\big[\delta_2'(Y+1)(Y+\frac{1}{2})Y\otimes1+\big(Y+1)(Y+\frac{1}{2})Y\otimes\delta_2']
=(\frac{-i\hbar}{2})^2\frac{1}{12}\big[\delta_2'(Y+1)(Y+\frac{1}{2})Y\otimes1
+\big(Y+1)(Y+\frac{1}{2})Y\otimes\delta_2']$.

\item $R_{0,0,0,0,3;0,0,3,0,0}=(\frac{-i\hbar}{2})^3\frac{i}{6\hbar}\big[\delta_2'\otimes
(Y+1)(Y+\frac{1}{2})Y+1\otimes \delta_2'(Y+1)(Y+\frac{1}{2})Y\big]
=(\frac{-i\hbar}{2})^2\frac{1}{12}\big[\delta_2'\otimes
(Y+1)(Y+\frac{1}{2})Y+1\otimes
\delta_2'(Y+1)(Y+\frac{1}{2})Y\big]$.

\item $R_{2,0,\cdots,0;0,0,2,0,0}=(\frac{-i\hbar}{2})^2\frac{1}{2}X^2\otimes
(Y+\frac{1}{2})Y$.

\item $R_{0,0,2,0,0;2,0,\cdots,0,0}=(\frac{-i\hbar}{2})^2\frac{1}{2}\Big((Y+\frac{1}{2})Y\otimes X^2
+\delta_2'(Y+\frac{1}{2})Y\otimes Y+
2\delta_1(Y+\frac{1}{2})Y\otimes
XY+\delta_2'(Y+\frac{1}{2})Y\otimes X+\delta_1^2(Y+\frac{1}{2})Y
\otimes Y^2+\frac{1}{2}\delta_1^2(Y+\frac{1}{2})Y\otimes
Y-\delta_2'Y(Y+\frac{1}{2})\otimes Y\Big)$.

\item $R_{1,0,1,0,0;1,0,1,0,0}=-(\frac{-i\hbar}{2})^2(X(Y+\frac{1}{2})\otimes X(Y+\frac{1}{2})
+\delta_1X(Y+\frac{1}{2})\otimes Y(Y+\frac{1}{2}))$
\end{itemize}

Taking the sum of all the above terms, we have that the $\hbar^2$
component of $R$ is equal to
\begin{eqnarray*}
&&(\frac{-i\hbar}{2})^2\Big(\frac{1}{2}\delta_2'
(Y+\frac{1}{2})Y\otimes Y+\frac{1}{6}\delta_2'\otimes
(Y+1)(Y+\frac{1}{2})Y+\frac{1}{2}X^2\otimes (Y+\frac{1}{2})Y\\
&&+\frac{1}{2}(Y+\frac{1}{2})Y\otimes X^2-X(Y+\frac{1}{2})\otimes
X(Y+\frac{1}{2})-\delta_1X(Y+\frac{1}{2})\otimes
Y(Y+\frac{1}{2})\\
&&+\delta_1(Y+\frac{1}{2})Y\otimes X(Y+\frac{1}{2})
+\frac{1}{2}\delta_1^2(Y+\frac{1}{2})Y\otimes
Y^2+\frac{1}{4}\delta_1^2(Y+\frac{1}{2})Y\otimes Y)\Big)
\end{eqnarray*}
\begin{eqnarray*}
&=&(\frac{-i\hbar}{2})^2\Big(\frac{1}{2}S(X)^2\otimes
(Y+\frac{1}{2})Y+S(X)(Y+\frac{1}{2})\otimes
X(Y+\frac{1}{2})+\frac{1}{2}Y(Y+\frac{1}{2})\otimes X^2\\
&&+\frac{1}{2}\delta_2' (Y+\frac{1}{2})Y\otimes
Y+\frac{1}{6}\delta_2'\otimes
(Y+1)(Y+\frac{1}{2})Y+\frac{1}{2}\delta_2'Y\otimes(Y+\frac{1}{2})Y\Big)
\end{eqnarray*}
\begin{remark}
The expression of $R_2$ agrees with the computation in
\cite{bty}[Prop. 6.1] up to a term
$\frac{1}{12}(\frac{-i\hbar}{2})^2\delta_2'\otimes Y$. We notice
that $\frac{1}{12}(\frac{-i\hbar}{2})^2\delta_2'\otimes Y$ is
Hochschild closed. For the purpose of \cite{bty}[Prop. 6.1], a
change of a closed Hochschild 2-cochain on $\hbar^2$ component
does not change the answer. This explains the difference.
\end{remark}

We remark that when $\delta_2'$ is an inner derivation, we can
replace $\delta_2'$ by $[\Omega,\cdot]$ in the expression of $R$.
But it turns out that the above computed $\hbar^2$ term $R_2$ does
not agree with $RC_2$ defined by Connes and Moscovici
\cite{CM03-2}, i.e. $R_2|_{\delta_2'[\Omega,\cdot]}-RC_2=$
\[
\Omega Y\otimes Y(2Y+1)+\Omega Y(2Y+1)\otimes
Y+\frac{1}{3}\Omega\otimes (Y+1)(2Y+1)Y-\frac{1}{3}\otimes \Omega
(Y+1)(2Y+1)Y.
\]

We do not know the explicit relation between these two universal
deformation formulas $R_2|_{\delta_2'=[\Omega, \cdot]}$ and
$RC_2$, but only a heuristic and geometric explanation to their
difference. The geometric constructions in \cite{bty} and in this
article are not exactly same. In \cite{bty}, we used a projective
structure to redefine a symplectic connection and a Fedosov
connection $D'$ on the Weyl algebra bundle $\calw$, and therefore
a symplectic diffeomorphism preserving this new Fedosov connection
naturally lifts to an endomorphism of the quantum algebra
$\calw_{D'}$. In this paper, we do not change the symplectic
connection because of lack of data, but change the quantization
process of a symplectic diffeomorphism by introducing sections
like $U_\alpha, U_\beta,\cdots$. Furthermore, we notice that the
difference $R_2|_{\delta'_2=[\Omega,\cdot]}-RC_2$ is actually a
Hochschild coboundary of a 1-Hochschild cochain
$-1/3\Omega(Y+1)(2Y+1)Y$. This suggests if we define an
isomorphism $I=1+1/3\hbar^2 \Omega(Y+1)(2Y+1)Y$ on
$C^\infty_c(\reals\times \reals^+)\rtimes \Gamma[[\hbar]]$, then
$I^{-1}(m(R(I(a)\otimes I(b))))=m+\hbar RC_1+\hbar^2
RC_2+o(\hbar^2)$. In general, we expect that if $\calh_1$ acts on
$\cala$ with a projective structure $\Omega$, there is an
isomorphism $I$ on $\cala[[\hbar]]$ which can be expressed using
elements in $\calh_1[[\hbar]]$ and the projective structure
$\Omega$ such that
\[
I^{-1}(m(R(I(a)\otimes I(b))))=m(RC(a\otimes b)).
\]

\section{Appendix : Associativity of the Eholzer Product}

In this appendix we study associativity of the Eholzer product,
which was used in Connes and Moscovici's approach \cite{CM03-2} to
obtain the general associativity at the Hopf algebra level. This
associativity theorem was first proved by Cohen, Manin, and Zagier
in \cite{CMZ}. In the first part of this appendix, we give a new
proof of the associativity using the method developed by the
second author \cite{yao}. In the second part, we study an
important combinatorial identity used by Cohen, Manin, and Zagier
in \cite{CMZ}. This interesting identity was obtained by Zagier
\cite{Z}, but its complete proof is missing in literature. We
prove this identity in the special case corresponding to the
Eholzer product.

\subsection{Proof of Associativity}

First we follow the argument developed in \cite{yao} according to
which the associativity of the product
\[f\ast g=\sum_{n=0}^\infty [f, g]_n \hbar^n\]
is equivalent to prove the identity (with the notation $
X_n=\prod_{i=0}^{n-1} (X+i)$)
\[
\sum_{r=0} {n-r \choose p} \frac{A_{n-r}(2k+2l+2r,  2m) A_{r}(2k,
2l) }{(2k+2l+2r)_{n-p-r} (2m)_p  (2k)_r  } =\sum_{s=0} {n-s
\choose n-p} \frac{A_{n-s}(2k,  2l+2m+2s) A_{s}(2l, 2m)
}{(2k)_{n-p}(2l+2m+2s)_{p-s} (2m)_s},
\]
for $p=0, 1, \dots,  n$ and for
\[
A_n(2k, 2l)= \frac{1}{n!}(2k)_n (2l)_n,
\]
i.e.
\begin{eqnarray}\label{ident}
& & \sum_{r=0}  {n-r \choose p}
\frac{\displaystyle\frac{1}{r!}(2k)_r(2l)_r}{(2k)_r}
\frac{\displaystyle\frac{1}{(n-r)!}(2k+2l+2r)_{n-r}(2m)_{n-r}}
{(2k+2l+2r)_{n-p-r}(2m)_p}\cr  & = &   \sum_{s=0}  {n-s \choose
n-p} \frac{\displaystyle\frac{1}{s!}(2l)_s(2m)_s}{(2m)_s}
\frac{\displaystyle\frac{1}{(n-s)!}(2k)_{n-s}(2l+2m+2s)_{n-s}}
{(2l+2m+2s)_{p-s}(2k)_{n-p}} .
\end{eqnarray}

Our proof is based on manipulation of combinatorial identities. We
have, for the left hand side,
\begin{eqnarray}\label{44}
& & \sum_{r=0}  {n-r \choose p}
\frac{\displaystyle\frac{1}{r!}(2k)_r(2l)_r}{(2k)_r}
\frac{\displaystyle\frac{1}{(n-r)!}(2k+2l+2r)_{n-r}(2m)_{n-r}}
{(2k+2l+2r)_{n-p-r}(2m)_p}\cr & = &  \sum_{r=0}
\frac{(n-r)!}{p!(n-r-p)!} \frac{1}{r!}\frac{(2k)_r(2l)_r}{(2k)_r}
\frac{1}{(n-r)!}
\frac{(2k+2l+2r)_{n-r}(2m)_{n-r}}{(2k+2l+2r)_{n-p-r}(2m)_p}\cr & =
& \sum_{r=0}
\frac{(2l)_r}{r!}\frac{(2k+2l+2r)_{n-r}}{(2k+2l+2r)_{n-p-r}p!}
 \frac{(2m)_{n-r}}{(2m)_p (n-r-p)!}\cr & = &
\sum_{r=0} {2l+r-1 \choose r} {2k+2l+n+r-1 \choose p}  {2m+n-r-1
\choose n-p-r}.
\end{eqnarray}

Once $n, p$  is fixed,  what to be verified is an identity about
polynomials in $2k, 2l, 2m$. When $2l$ is a negative integer, we
can use the following two combinatorial relations:
\begin{eqnarray*}
{X+n-1\choose n} =(-1)^n{-X\choose n}, n>0,\ \ \ \text{and\ }\ \ \
\sum_{i}{X\choose i}{Y\choose n-i}={X+Y\choose n},
\end{eqnarray*}
to get
\begin{eqnarray*}
{2l+r-1 \choose r} & = & (-1)^r {-2l \choose r}, \cr {2k+2l+n+r-1
\choose p} &=& \sum_u {2k+2l+n-1 \choose p+2l+u}{r \choose -2l-u},
\cr {2m+n-r-1 \choose n-p-r} & = & \sum_v {2l+2m+n-1 \choose
n-p+2l+v}{-2l-r \choose -2l-r-v}.
\end{eqnarray*}

Then (\ref{44}) becomes {\small
\begin{eqnarray*}
&= &  \sum_{r=0} (-1)^r {-2l \choose r}  \left[\sum_u {2k+2l+n-1
\choose p+2l+u}{r \choose -2l-u}\right]\left[\sum_v {2l+2m+n-1
\choose n-p+2l+v}{-2l-r \choose -2l-r-v}\right]\cr
&= & \sum_{u,
v} {2k+2l+n-1\choose p+2l+u} {2l+2m+n-1 \choose
n-p+2l+v}\left[\sum_r (-1)^r {-2l\choose r}{r\choose
-2l-u}{-2l-r\choose -2l-r-v}\right].
\end{eqnarray*}}

We then simplify the quantity inside the above brackets by
\begin{eqnarray*}
& & \sum_r (-1)^r {-2l\choose r}{r\choose -2l-u}{-2l-r\choose
-2l-r-v} \cr &=& \sum_r (-1)^r
\frac{(-2l)!}{r!(-2l-r)!}\frac{r!}{(-2l-u)!(r+2l+u)!}\frac{(-2l-r)!}{(
-2l-r-v)!v!}\cr &=&\frac{(-2l)!}{(-2l-u)!v!}\sum_r (-1)^r
\frac{1}{(r+2l+u)!(-2l-r-v)!}\cr &=&
\frac{(-2l)!}{(-2l-u)!v!}\frac{1}{(u-v)!}\sum_r (-1)^r
\frac{(u-v)!}{(r+2l+u)!(-2l-r-v)!}\cr &=&
\frac{(-2l)!}{(-2l-u)!v!}\frac{1}{(u-v)!}\sum_r (-1)^r {u-v\choose
-2l-r-v}\cr &=&
\frac{(-2l)!}{(-2l-u)!v!}\frac{1}{(u-v)!}(1-1)^{u-v}(-1)^{-2l-v}=\frac{(-2l)!}{(-2l-u)!
v!}(-1)^{-2l-v}\delta_{u, v},
\end{eqnarray*}
where $\delta_{x, y}$ is the Kronecker symbol ($\delta_{x, y}=1$
if $x=y$,  $=0$ if not). Finally we get,
\begin{eqnarray*}
& & \sum_{r=0}  {n-r \choose p}
\frac{\displaystyle\frac{1}{r!}(2k)_r(2l)_r}{(2k)_r}
\frac{\displaystyle\frac{1}{(n-r)!}(2k+2l+2r)_{n-r}(2m)_{n-r}}
{(2k+2l+2r)_{n-p-r}(2m)_p}\cr & = &  \sum_{u=v=-2l-t}
{2k+2l+n-1\choose p+2l+u} {2l+2m+n-1 \choose n-p+2l+v}
\frac{(-2l)!}{(-2l-u)! v!}(-1)^{-2l-v}\cr & = & \sum_t
{2k+2l+n-1\choose p-t} {2l+2m+n-1 \choose n-p-t}(-1)^t {-2l
\choose t} .
\end{eqnarray*}

On the right hand side of (\ref{ident}),  we have
\begin{eqnarray*}
& & \sum_{s=0}  {n-s \choose n-p}
\frac{\displaystyle\frac{1}{s!}(2l)_s(2m)_s}{(2m)_s}
\frac{\displaystyle\frac{1}{(n-s)!}(2k)_{n-s}(2l+2m+2s)_{n-s}}
{(2l+2m+2s)_{p-s}(2k)_{n-p}}\cr & = &\sum_{s=0}
\frac{(n-s)!}{(n-p)!(p-s)!}\frac{1}{s!}
\frac{(2l)_s(2m)_s}{(2m)_s}
\frac{1}{(n-s)!}\frac{(2k)_{n-s}(2l+2m+2s)_{n-s}}{(2k)_{n-p}(2l+2m+2s)_{p-s}}\cr
& = & \sum_{s=0} {2l+s-1\choose s}{2k+n-s-1\choose
p-s}{2l+2m+s+n-1\choose n-p}.
\end{eqnarray*}

By the same method as before, we compute the above quantity  as
follows,{ \small
\begin{eqnarray*} & =& \sum_{s=0} (-1)^s {-2l\choose s}
\left[\sum_v {2k+2l+n-1\choose p+2l+v}{-2l-s\choose -2l-s-v}
\right]\ \left[\sum_u {2l+2m+n-1\choose n-p+2l+u}{s\choose
-2l-u}\right]\cr & = & \sum_{u, v} {2k+2l+n-1\choose p+2l+v}
{2l+2m+n-1\choose n-p+2l+u}\left[\sum_{s-0} (-1)^s {-2l\choose
s}{-2l-s\choose -2l-s-v} {s\choose -2l-u}\right]\cr & =&
\sum_{u=v=-2l-t} {2k+2l+n-1\choose p+2l+v} {2l+2m+n-1\choose
n-p+2l+u} \frac{(-2l)!}{(-2l-u)!v!}(-1)^{-2l-v}\cr &=& \sum_t
{2k+2l+n-1\choose p-t}{2l+2m+n-1\choose n-p-t}(-1)^t{-2l\choose t}
\end{eqnarray*}}
which gives out the same quantity. We obtain then

\begin{proposition}
The Eholzer product is associative.
\end{proposition}
\begin{remark}
If we want to prove the identification of the coefficients of
every $d^{u}f d^{v}g d^{w}h$ in $(f\ast g)\ast h$ and $f\ast(g\ast
h)$, what we proved above is that for every triple of indices
$(l_1, l_2, l_3)$, the following combinatorial identity :

\begin{eqnarray*}
& & \displaystyle \sum_{r=0}^{l_1}\sum_{s=0}^{l_2}\sum_{t=0}^{l_3}
(-1)^{l_1+l_2-s} {E+l_1+s+l_3-t-1 \choose l_3-t}{E+r+s-1 \choose
s}f \cr & & \hspace{2cm}  {E+l_2+r+t-1\choose t}{E+r+s-1\choose
r}g\cr & & \hspace{3cm}
{E+l_1-r+l_2-s+l_3-1\choose l_1-r}{E+l_3+l_2-s-1\choose l_2-s}h\\
 \cr &=&  \sum_{r=0}^{l_1}\sum_{s=0}^{l_2}\sum_{t=0}^{l_3}
(-1)^{l_1+l_2-s} {E+l_1+s+l_3-t-1\choose l_3-t}{E+l_1+s-1\choose
s}f\cr & & \hspace{2cm} {E+l_2-s+t-1\choose t}{E+l_2+t+r-1\choose
r}g\cr & & \hspace{3cm} {E+l_1-r+l_2-s+l_3-1\choose
l_1-r}{E+t+l_2-s-1\choose l_2-s}h,
\end{eqnarray*}

\medskip
\noindent where $E$ is the Euler operator.
\end{remark}

\subsection{Zagier's identity}

In this part we turn to the original Cohen, Manin and Zagier's
proof \cite{CMZ} of the Eholzer product. Their proof relies on the
following combinatorial identity

\begin{equation}\label{cent}
\frac{(-4)^n}{\displaystyle{2x \choose n}}\sum_{r+s=n}
\frac{\displaystyle{y \choose r}{y-a \choose r}}{\displaystyle{2y
\choose r}}\frac{\displaystyle{z\choose s}{z+a\choose
s}}{\displaystyle{2z\choose s}}=\sum_{j\geq 0}{n \choose
2j}\frac{\displaystyle{-\frac{1}{2}\choose j}{a-\frac{1}{2}\choose
j}{-a-\frac{1}{2}\choose j}}{\displaystyle{x-\frac{1}{2}\choose
j}{y-\frac{1}{2}\choose j}{z-\frac{1}{2}\choose j}}
\end{equation}
where $n\geq 0$ and the variables $a, x, y, z$ satisfy
$x+y+z=n-1$.

In this part, we give a proof of this identity when
$a=\frac{1}{2}$. We start with some transformations. Our aim is to
eliminate the binomial coefficients in the denominator of both
sides. Using the following identity
\[
\displaystyle\frac{1}{\displaystyle{X\choose
r}}=\frac{\displaystyle{X-r\choose n-r}}{\displaystyle{X\choose
n}}\frac{r!(n-r)!}{n!},
\]
we rewrite the left hand side of (\ref{cent}) into the following
expression
\begin{eqnarray}\label{gauche}
\frac{(-4)^n}{\displaystyle{2x\choose n}} \sum_{r=0}^n
\frac{\displaystyle{y\choose r}{y-a\choose r }{2y-r\choose
n-r}}{\displaystyle{2y\choose n}} \frac{\displaystyle{z\choose
n-r}{z+a \choose n-r}{2z-n+r\choose r}}{\displaystyle{2z\choose
n}}\frac{(r!(n-r)!)^2}{(n!)^2}.\cr\hspace{-1cm} & &
\end{eqnarray}

And using the following identity
\[
{2Y\choose 2j}={Y-\frac{1}{2} \choose j}{Y\choose j}\frac{(j!)^2
4^j}{(2j)!},
\]
we rewrite the right hand side of (\ref{cent}) into the following
expression,
\begin{eqnarray}\label{droit}
&= & \sum_{j=0}^{\left[\frac{n}{2}\right]} {n\choose
2j}{-\frac{1}{2}\choose j}{a-\frac{1}{2}\choose
j}{-a-\frac{1}{2}\choose j} \frac{\displaystyle{x\choose j} 4^j
(j!)^2}{\displaystyle{2x\choose 2j
}(2j)!}\frac{\displaystyle{y\choose j} 4^j
(j!)^2}{\displaystyle{2y\choose 2j
}(2j)!}\frac{\displaystyle{z\choose j} 4^j
(j!)^2}{\displaystyle{2z\choose 2j }(2j)!}\cr & = &
\sum_{j=0}^{\left[\frac{n}{2}\right]} {n\choose
2j}{-\frac{1}{2}\choose j}{a-\frac{1}{2}\choose
j}{-a-\frac{1}{2}\choose j} \left[\frac{4^j
(j!)^2(n-2j)!}{n!}\right]^3 \cr & &\hspace{1cm}
\frac{\displaystyle{x\choose j}{2x-2j\choose
n-2j}}{\displaystyle{2x\choose n}}\frac{\displaystyle{y\choose
j}{2y-2j\choose n-2j}}{\displaystyle{2y\choose
n}}\frac{\displaystyle{z\choose j}{2z-2j\choose
n-2j}}{\displaystyle{2z\choose n}}.
\end{eqnarray}

By combining (\ref{gauche}) and (\ref{droit}) we can then multiply
both sides of (\ref{cent}) by the common denominator ${2x\choose
n}{2y\choose n}{2z\choose n} $. We obtain a polynomial $P_n(y, z,
a)$ on the left hand side of (\ref{cent}),
\begin{eqnarray*}
P_n(y, z, a)&:=& (-4)^n \sum_{r=0}^n {y\choose r}{y-a\choose r
}{2y-r\choose n-r}{z\choose n-r}{z+a \choose n-r}{2z-n+r\choose r}
(r!(n-r)!)^2, \end{eqnarray*}
and a polynomial $Q_n(y, z, a)$ on
the right hand side of (\ref{cent}), (we replace $x$ by $n-y-z-1$)
\begin{eqnarray*}\label{Qndef}
    Q_n(y, z, a)&:=&\sum_{j=0}^{\left[\frac{n}{2}\right]}
    \frac{((n-2j)!)^2 (j!)^6 2^{6j}}{(2j)!}{-\frac{1}{2}\choose
    j}{a-\frac{1}{2}\choose j}{-a-\frac{1}{2}\choose j} \cr & &
    \hspace{0.8cm}{n-y-z-1\choose j }{2n-2y-2z-2-2j\choose n-2j}\cr & &
    \hspace{1.6cm}{y\choose j}{2y-2j
    \choose
    n-2j}{z\choose j}{2z-2j\choose n-2j}. \cr & &
\end{eqnarray*}

In summary, identity (\ref{cent}) is equivalent to the identity $
`` P_n(y,z,a)=Q_n(y,z,a) " $ for $n\geq 0$. In order to have the
associativity of the Eholzer product, we shall prove this identity
when $a=\frac{1}{2}$. Explicitly, we want to prove the following
identity,
\begin{equation}\label{k=01}
\begin{split}
&(-1)^n \sum_{r=0}^n {2y\choose 2r} {2y-r\choose n-r} {2z+1\choose
2(n-r)}{2z-n+r\choose
r}\frac{(2r)!(2(n-r))!}{(n!)^2}\\
&={2n-2y-2z-2\choose n}{2y \choose n}{2z\choose n }.
\end{split}
\end{equation}

\subsubsection{Simplification}
We apply the following identity on the left hand side of
(\ref{k=01})
\begin{eqnarray*}
\frac{\displaystyle{2y\choose 2r}{2y-r\choose n-r}(2r)!}{n!} & =&
{2y\choose n}{2y-r\choose r}\frac{r!}{(n-r)!},\cr
\frac{\displaystyle{2z+1\choose 2(n-r)}{2z-n+r\choose
r}(2(n-r))!}{n!}&=&{2z\choose n}\left[{2z-n+r\choose
n-r}+2{2z-n+r\choose n-r-1}\right] \frac{(n-r)!}{r!}.
\end{eqnarray*}

Taking the quotients on both sides of (\ref{k=01}) by ${2y\choose
n}{2z\choose n}$, we have the following equivalent identity,
\[
\hspace{-2cm} (-1)^n \sum_{r=0}^n {2y-r\choose r}
\left[{2z-n+r\choose n-r}+2{2z-n+r\choose n-r-1}\right]=
{2n-2y-2z-2\choose n}.
\]
And we can rewrite the above identity as
\begin{equation}\label{1/2}
\hspace{-1cm}\sum_{r=0}^n {2y-r\choose r} \left[{2z-n+r\choose
n-r}+2{2z-n+r\choose n-r-1}\right]={2y+2z-n+1\choose n},
\end{equation}
using ${2n-2y-2z-2\choose n}=(-1)^n {2y+2z-n+1\choose n}$.

\subsubsection{The sums $S_0(n;A, B)$ and $S(n;X)$}

We consider the sum
\[
S_0(n;A, B)=\sum_{k=0}^{n} {k+A \choose n-k} {n-k+B\choose k}.
\]
We have $S_0(0;A, B)=1$, and by $ {X+1\choose n}={X\choose
n}+{X\choose n-1}$,
\begin{eqnarray*}
\displaystyle S_0(n+1;A, B)& =& \sum_{k=0}^{n+1} {k+A \choose
n+1-k} {n+1-k+B\choose k} \cr & = &
 S_0(n;A-1, B+1)+S_0(n+1;A-1, B).
\end{eqnarray*}
In the same way, $S_0(n+1;A, B)=S_0(n;A+1, B-1)+S_0(n+1;A, B-1)$.

On the other hand, we define now
\[
S(n;X)=\sum_{p=0}^{\left[\frac{n}{2}\right]} {X+n-1-2p\choose
n-2p}.
\]
It is clear that $S_0(0;X)=1$ and $S(n+1;X)=S(n+1;X-1)+S(n;X)$ for
the same reason as above.

Another useful relation is that
\begin{equation}\label{n, n-1}
 S(n, X-1)+2S(n-1, X)={X+n\choose n},
\end{equation}
because
\begin{eqnarray*}
S(n, X-1)+2S(n-1, X)&=& S(n;X)+S(n-1;X)\cr
&=&\sum_{p=0}^{\left[\frac{n}{2}\right]}{X+n-1-2p\choose n-2p}+
\sum_{p=0}^{\left[\frac{n-1}{2}\right]}{X+n-2-2p\choose n-1-2p}\cr
&=& \sum_{i=0}^n (-1)^i {-X\choose i}= (-1)^n {-X-1\choose
n}={X+n\choose n}.
\end{eqnarray*}

We prove the following lemma:
\begin{lemma}\label{s0=s}
$S_0(n;A, B)=S(n;A+B)$.
\end{lemma}

\noindent {\bf Proof.} We prove the lemma by induction on $n$.
When $n=0$,  we have obviously equality. If the claim is valid for
$0,\, 1,\, \cdots,\, n$, then by the above two induction
relations, we have
\begin{eqnarray*}
\hspace{-1cm}& & S_0(n+1;A,\, B)-S(n+1;A+B)\cr \hspace{-1cm}&=&
(S_0(n+1;A,\,B-1)+S_0(n;A+1,\, B-1))
-(S(n+1;A+B-1)+S(n;A+B))\cr\hspace{-1cm} &=& (S_0(n+1;A,\,
B-1)-S(n+1;A+B-1)) +(S_0(n;A+1,\, B-1)-S(n;A+B))\cr\hspace{-1cm}
&=& S_0(n+1;A,\, B-1)-S(n+1;A+B-1).
\end{eqnarray*}

With the same method, we can prove that this difference also
equals $S_0(n+1;A-1,\, B)-S(n+1;A+B-1)$. We hence conclude that
the difference $S_0(n+1;A,B)-S(n+1;A+B)$ has the same value for
all pairs $(A,\, B)\in {\mathbb N}^2$. But by definition we can
observe
\[S(n+1;0,\, 0)-S(n+1;0)=\left\{\begin{array}{cl}
                  0-0=0  &  \text{if}\, \, n\, \,  \text{odd}, \\
                  1-1=0  &  \text{if}\, \, n\, \,  \text{even}.
                     \end{array}\right.\]
We conclude that the identity holds for all $n$ and all pairs
$(A,\, B)\in {\mathbb N}^2$. But as what we want to prove is a
polynomial identity in $A$ and $B$(for fixed $n$), the identities
for all natural numbers implies its correctness for arbitrary
$(A,\, B)$.\ $\Box$

\subsubsection{Resummation}

\begin{theorem}
The identity (\ref{1/2}) is valid.
\end{theorem}

\noindent{\bf Proof.} In fact, by using (\ref{n, n-1}), we have
\begin{eqnarray*}
& &\sum_{r=0}^n {2y-r\choose r} \left[{2z-n+r\choose
n-r}+2{2z-n+r\choose n-r-1}\right]\cr & = & [S_0(n;2z-n,
2y-n)+2S_0(n-1;2z-n, 2y-n+1)]\cr &
=&[S(n;2y+2z-2n)+2S(n-1;2y+2z-2n+1)]\cr &=& { 2y+2z-n+1 \choose
n}.\ \ \Box
\end{eqnarray*}

\vspace{2mm}

{\small \noindent{Xiang Tang}, Department of Mathematics,
Washington University, St. Louis, MO, 63139, U.S.A.,\\
Email:xtang@math.wustl.edu.

\vspace{2mm}

\noindent{Yijun Yao}, Projet AO, IMJ, 175 Rue du Chevalret, 75013
Paris, France, email:yao@math.jussieu.fr.}

\begin{thebibliography}{99}
\bibitem{bty}Bieliavsky, P., Tang, X., and Yao, Y., Rankin-Cohen
brackets and formal deformation, to appear in {\em Adv.
Mathematics}.

\bibitem{bgnt}Brussler, A., Gorokhovsky, A., Nest, R., and Tsygan,
B., Deformation quantization of gerbes, to appear in {\em Adv.
Mathematics}.

\bibitem{CM03-1}Connes, A., and Moscovici, H.,  Modular , Hecke algebras and their Hopf
symmetry, {\em Mosc. Math. J.} 4 (2004), no. 1, 67--109.

\bibitem{CM03-2}Connes, A., and Moscovici, H., Rankin-Cohen brackets and the Hopf algebra
 of transverse geometry,  {\em Mosc. Math. J.}  4  (2004),  no. 1, 111--130.

\bibitem{CMZ}Cohen, P., Manin, Y., and Zagier, D.,
Automorphic pseudodifferential operators,  {\em Algebraic aspects
of integrable systems},  17--47, Progr. Nonlinear Differential
Equations Appl., 26, Birkha\"user Boston, Boston, MA, 1997.


\bibitem{fe}Fedosov, B., A simple geometrical construction of deformation quantization,
{\em J. Differential Geom.}  40  (1994),  no. 2, 213--238.

\bibitem{gi-zh}Giaquinto, A., and Zhang, J., Bialgebra actions, twists, and
universal deformation formulas,  {\em J. Pure Appl. Algebra}  128
(1998),  no. 2, 133--151.

\bibitem{yao} Yao, Y., Rankin-Cohen Deformations and Representation
Theory, {\em arxiv:math.0708.1528}.

\bibitem{Z94} Zagier,  D.,  Modular forms and differential operators. K. G.
Ramanathan memorial issue,   {\em Proc. Indian Acad. Sci. Math.
Sci.} 104 (1994),   no. 1,  57--75.

\bibitem{Z01} Zagier,  D.,  Formes modulaires et Op\'erateurs diff\'erentiels,
\textit{2001-2002 Course at Coll\`ege de France}.

\bibitem{Z} Zagier,  D.,  Some combinatorial identities occuring in
the theory of modular forms,  in preparation.
\end{thebibliography}
\end{document}